\documentclass[9pt]{article}
\usepackage{amsthm,amsfonts,amsmath,amssymb,amscd}
\usepackage[arc,all,knot,poly]{xy}
\usepackage{bm}
\usepackage[dvips]{color}
\theoremstyle{definition}
\newtheorem{thm}{Theorem} [section]
\newtheorem{cor}[thm]{Corollary}
\newtheorem{lem}[thm]{Lemma}

\newtheorem{rem}[thm]{Remark}
\newtheorem{defi}[thm]{Definition}
\newtheorem{exm}[thm]{Example}
\author{Hitoshi Yamanaka}
\title{\textbf{Weights of Markov traces for \linebreak Alexander polynomials of mixed links }}
\date{}
\begin{document}
\maketitle
\begin{abstract}
Using the Fourier expansion of Markov traces for Ariki-Koike algebras over $\mathbb{Q}(q,u_{1},\dots,u_{e})$, we give a direct definition of the Alexander polynomials
for mixed links. We observe that under the corresponding specialization of a Markov parameter, the Fourier coefficients of Markov traces take quite simple form.

As a consequence, we show that the Alexander polynomial of a mixed link is essentially equal to the Alexander polynomial of the link obtained by resolving the
twisted parts. 


 
\end{abstract}
\section{Introduction}
In [L1] Lambropoulou initiated the mixed link theory which is the link theory in the solid torus. A mixed link is an embedding of a disjoint union of finitely many circles into the solid torus $Y$. Since the 3-sphere has the canonical
genus 1 Heegaard decomposition, it can be considered as an embedding into the complement $S^{3}\setminus \text{Int}\ Y$.
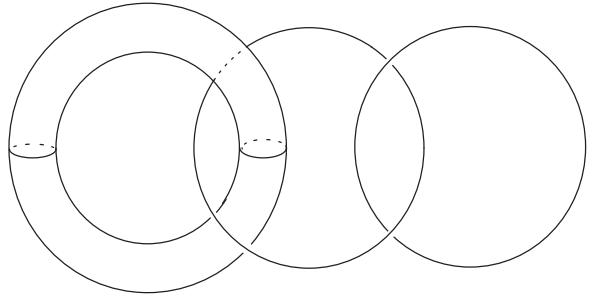
\begin{center}
\unitlength 0.1in
\begin{picture}( 30.1600, 15.1400)(  4.0400,-19.2000)
%
{\color[named]{Black}{%
\special{pn 8}%
\special{ar 1130 1164 726 758  0.8137996 6.2831853}%
\special{ar 1130 1164 726 758  0.0000000 0.7259561}%
}}%
%
{\color[named]{Black}{%
\special{pn 8}%
\special{ar 1130 1164 480 502  0.8188646 6.2831853}%
\special{ar 1130 1164 480 502  0.0000000 0.6168629}%
}}%
%
{\color[named]{Black}{%
\special{pn 8}%
\special{ar 528 1170 122 46  6.1517166 6.2831853}%
\special{ar 528 1170 122 46  0.0000000 3.1415927}%
}}%
%
{\color[named]{Black}{%
\special{pn 8}%
\special{ar 1734 1164 122 52  6.2831853 6.2831853}%
\special{ar 1734 1164 122 52  0.0000000 3.1415927}%
}}%
%
{\color[named]{Black}{%
\special{pn 8}%
\special{ar 528 1176 122 32  3.1415927 3.2994874}%
\special{ar 528 1176 122 32  3.7731716 3.9310663}%
\special{ar 528 1176 122 32  4.4047505 4.5626453}%
\special{ar 528 1176 122 32  5.0363295 5.1942242}%
\special{ar 528 1176 122 32  5.6679084 5.8258032}%
}}%
%
{\color[named]{Black}{%
\special{pn 8}%
\special{ar 1740 1164 116 44  3.0152845 3.1662279}%
\special{ar 1740 1164 116 44  3.6190581 3.7700015}%
\special{ar 1740 1164 116 44  4.2228317 4.3737751}%
\special{ar 1740 1164 116 44  4.8266053 4.9775487}%
\special{ar 1740 1164 116 44  5.4303789 5.5813223}%
\special{ar 1740 1164 116 44  6.0341525 6.1850959}%
}}%
%
{\color[named]{Black}{%
\special{pn 8}%
\special{ar 1974 1164 604 630  4.1375807 5.2887642}%
}}%
%
{\color[named]{Black}{%
\special{pn 8}%
\special{ar 1974 1164 602 628  3.7402628 3.7597908}%
\special{ar 1974 1164 602 628  3.8183750 3.8379031}%
\special{ar 1974 1164 602 628  3.8964873 3.9160154}%
\special{ar 1974 1164 602 628  3.9745996 3.9941277}%
\special{ar 1974 1164 602 628  4.0527119 4.0722400}%
}}%
%
{\color[named]{Black}{%
\special{pn 8}%
\special{ar 1974 1164 602 630  5.5274344 6.2831853}%
\special{ar 1974 1164 602 630  0.0000000 3.7402628}%
}}%
%
{\color[named]{Black}{%
\special{pn 8}%
\special{ar 2818 1158 604 630  2.3662613 6.2831853}%
\special{ar 2818 1158 604 630  0.0000000 2.2927503}%
}}%
%
{\color[named]{Black}{%
\special{pn 8}%
\special{ar 1974 1164 602 630  5.2842267 5.4482993}%
}}%
%
{\color[named]{Black}{%
\special{pn 8}%
\special{ar 1130 1158 488 508  0.5700854 0.7406949}%
}}%
\end{picture}%

\end{center}
We say that two mixed links are equivalent if they are joined by an ambient isotopy of $S^{3}\setminus \text{Int}\ Y$.

Lambropoulou showed that there is a variant of the usual braid theory, that is, every mixed link is equivalent to
the closure of a mixed braid and one can consider analogues of the Markov Moves. More precisely, a mixed braid with $n$-strands is an $n$-tuple $(p_{1},\dots,p_{n})$ of embeddings of the closed interval $[0,1]$ into the closure $\text{cl}(I^{3}\setminus Cyn)$,
where 
\begin{center}
$\displaystyle Cyn:=\bigg\{(x,y,z)\in I^{3} \bigg| \bigg( x-\frac{1}{4}\bigg)^{2}+\bigg( y-\frac{1}{2(n+1)}\bigg)^{2}\leq \bigg(\frac{1}{4(n+1)}\bigg)^{2}\bigg\}$,
\end{center}
such that
\begin{itemize}
\item the curves $p_{i}(I)$ and $p_{j}(I)$ do not intersect if $i\not= j$,
\item there exists a permutation $\sigma $ of $\{ 1,\dots,n\}$ such that 
\begin{center}
$\displaystyle p_{i}(0)=\bigg(\frac{1}{2},\frac{i}{n+1},0\bigg),\ p_{i}(1)=\bigg(\frac{1}{2},\frac{\sigma(i)}{n+1},1\bigg)$
\end{center}
for all $i\in \{ 1,\dots, n\}$,
\item the point $p_{i}(t)$ lies in the interior of $\text{cl}(I^{3}\setminus Cyn)$ if $t\in (0,1)$,
\item the 3rd coordinate of $p_{i}(t)$ is increasing with respect to $t$ for all $i\in \{ 1,\dots, n\}$.
\end{itemize}
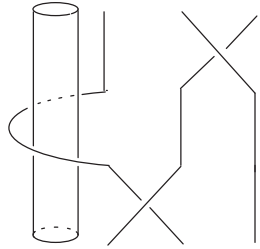
\begin{center}
\unitlength 0.1in
\begin{picture}( 12.9900, 12.7000)( 12.1100,-22.1000)
%
{\color[named]{Black}{%
\special{pn 8}%
\special{ar 1456 974 118 34  0.8934099 6.2831853}%
\special{ar 1456 974 118 34  0.0000000 0.7853982}%
}}%
%
{\color[named]{Black}{%
\special{pn 8}%
\special{ar 1456 2156 118 40  6.1049178 6.2831853}%
\special{ar 1456 2156 118 40  0.0000000 3.5204240}%
}}%
%
{\color[named]{Black}{%
\special{pn 8}%
\special{ar 1456 2156 118 40  2.8250234 2.9778896}%
\special{ar 1456 2156 118 40  3.4364884 3.5893546}%
\special{ar 1456 2156 118 40  4.0479533 4.2008196}%
\special{ar 1456 2156 118 40  4.6594183 4.8122845}%
\special{ar 1456 2156 118 40  5.2708833 5.4237495}%
\special{ar 1456 2156 118 40  5.8823482 6.0352145}%
}}%
%
{\color[named]{Black}{%
\special{pn 8}%
\special{pa 1730 2210}%
\special{pa 2112 1798}%
\special{fp}%
}}%
%
{\color[named]{Black}{%
\special{pn 8}%
\special{pa 2124 2204}%
\special{pa 1950 2014}%
\special{fp}%
}}%
%
{\color[named]{Black}{%
\special{pn 8}%
\special{pa 1736 1798}%
\special{pa 1904 1974}%
\special{fp}%
}}%
%
{\color[named]{Black}{%
\special{pn 8}%
\special{pa 2500 2196}%
\special{pa 2500 1792}%
\special{fp}%
}}%
%
{\color[named]{Black}{%
\special{pn 8}%
\special{ar 1918 1596 706 204  1.8311145 3.6687386}%
}}%
%
{\color[named]{Black}{%
\special{pn 8}%
\special{ar 1918 1596 714 196  4.2429049 4.4259287}%
}}%
%
{\color[named]{Black}{%
\special{pn 8}%
\special{ar 1918 1590 720 190  3.8352724 3.8617041}%
\special{ar 1918 1590 720 190  3.9409992 3.9674309}%
\special{ar 1918 1590 720 190  4.0467261 4.0731578}%
\special{ar 1918 1590 720 190  4.1524530 4.1788847}%
}}%
%
{\color[named]{Black}{%
\special{pn 8}%
\special{pa 1340 988}%
\special{pa 1340 1690}%
\special{fp}%
}}%
%
{\color[named]{Black}{%
\special{pn 8}%
\special{pa 1340 1744}%
\special{pa 1340 2156}%
\special{fp}%
}}%
%
{\color[named]{Black}{%
\special{pn 8}%
\special{pa 1574 982}%
\special{pa 1574 1752}%
\special{fp}%
}}%
%
{\color[named]{Black}{%
\special{pn 8}%
\special{pa 1574 2156}%
\special{pa 1574 1798}%
\special{fp}%
}}%
%
{\color[named]{Black}{%
\special{pn 8}%
\special{pa 1722 1400}%
\special{pa 1722 1400}%
\special{fp}%
}}%
%
{\color[named]{Black}{%
\special{pn 8}%
\special{pa 1710 1400}%
\special{pa 1710 988}%
\special{fp}%
}}%
%
{\color[named]{Black}{%
\special{pn 8}%
\special{pa 2500 1414}%
\special{pa 2500 1826}%
\special{fp}%
}}%
%
{\color[named]{Black}{%
\special{pn 8}%
\special{pa 2500 1414}%
\special{pa 2118 988}%
\special{fp}%
}}%
%
{\color[named]{Black}{%
\special{pn 8}%
\special{pa 2112 1792}%
\special{pa 2112 1394}%
\special{fp}%
}}%
%
{\color[named]{Black}{%
\special{pn 8}%
\special{pa 2110 1390}%
\special{pa 2286 1222}%
\special{fp}%
}}%
%
{\color[named]{Black}{%
\special{pn 8}%
\special{pa 2510 1010}%
\special{pa 2350 1180}%
\special{fp}%
}}%
\end{picture}%

\end{center}
Two mixed braids with $n$-strands are said to be equivalent if they are joined by an ambient isotopy of $\text{cl}(I^{3}\setminus Cyn)$ which fixes the boundary.

Then the set of equivalence classes is in bijection with the affine braid group $B_{\text{aff},n}$ with generators 
\begin{center}
$t_{0},t_{1},\cdots,t_{n-1}$
\end{center}
satisfying fundamental relations  
\begin{center}
$t_{0}t_{1}t_{0}t_{1}=t_{1}t_{0}t_{1}t_{0},$\\
$t_{i}t_{j}=t_{j}t_{i}\ \ \ \ \ (|i-j|\geq 2),$\\
$t_{i}t_{i+1}t_{i}=t_{i+1}t_{i}t_{i+1}\ \ \ \ \ (1\leq i\leq n-2)$.
\end{center}
It is known that $B_{\text{aff},n}$ is the semi-direct product of the braid group $B_{n}=\left<t_{1},\dots,t_{n-1}\right>$ and the free subgroup $P_{n}$
generated by $t_{0}',t_{1}'\dots,t_{n-1}'$ where 
\begin{center}
$t_{i}'=t_{i}\cdots t_{1}t_{0}t_{1}^{-1}\cdots t_{i}^{-1}\ \ \ \ \ (i=0,1,\dots,n-1)$.
\end{center}
In particular, the braid group $B_{n}$ is embedded in $B_{\text{aff},n}$ as a subgroup.

Moreover, in this setting we have the following analogues of the Markov moves:
\begin{itemize}
\item[(1)] $\alpha \longleftrightarrow \beta \alpha \beta ^{-1}\ \ (\alpha ,\beta \in B_{\text{aff},n})$,
\item[(2)] $\alpha \longleftrightarrow \alpha t_{n}^{\pm 1}\ \ (\alpha \in B_{\text{aff},n})$.
\end{itemize}

One of the main interest in this paper is to construct an analogue of the Alexander polynomial explicitly.
\bigskip

In the usual link theory, Jones [J] discovered a way to construct the HOMFLYPT polynomial using the Markov traces of Iwahori-Hecke algebras
of type $A$. He gave two methods to derive the Alexander polynomials.
The first one is to use the Skein relation for HOMFLYPT polynomials and the second is
to use the ``Fourier expansion" of the Markov traces of Iwahori-Hecke algebras of type $A$,
namely, the expression as the linear combination of irreducible characters. Note that the second method is more direct than the first one.

As for the mixed links, it is possible to give an analogue of the Alexander polynomial using the analogue of HOMFLYPT polynomials and their
Skein relations given in [L2]. We point out that one can also define the Alexander polynomial of a mixed link directly following Jones's second argument.
For this purpose a result of Geck-Iancu-Malle [GIM] is quite helpful. In [GIM], they determined the rational polynomials appearing 
in the coefficients in the expression of the Markov traces as the linear combination of the irreducible characters of the Ariki-Koike algebra of type $G(e,1,n)$.

Then our second observation is that when we consider the specialization of parameter, the Fourier coefficients take quite simple form. As a consequence,
we show that the Alexander polynomial of a mixed link is essentially the same as the Alexander polynomial of the link obtained by resolving 
the twisting parts. 
\bigskip

This paper is organized as follows. In Section 2, we recall the definition of the Ariki-Koike algebras and its irreducible ordinary representations.
In Section 3, we recall the definition of the Markov traces for Ariki-Koike algebras and its Fourier expansion due to Geck-Iancu-Malle.
In Section 4, we collect some  fundamental definitions and results in Lambropoulou's mixed link theory. 
In Section 5, we propose a definition of Alexander polynomial for a mixed link. In Section 6, we calculate the Fourier coefficient of Markov traces 
and prove a relation between the Alexander polynomials for mixed links and the one for the usual links.

\bigskip
\noindent
{\bf Acknowledgment.} The work was supported by Grant-in-Aid for JPSP Fellows $23\cdot 9693$.

\section{Representations of Ariki-Koike algebras}
In this section we recall the definition of the Ariki-Koike algebra and its ordinary finite dimensional irreducible representations.
For the details we refer the original paper [AK]. 

Let $e,n$ be two positive integers. We denote by $\Bbbk$  the field of rational functions over $\mathbb{Q}$ with $(e+1)$-indeterminates $q,u_{1},\cdots,u_{e}$.
\begin{defi}
The \textbf{Ariki-Koike algebra} of type $G(e,1,n)$ is the associative $\Bbbk$-algebra $H_{e,n}$ with generators $T_{0},T_{1},\cdots,T_{n-1}$
satisfying fundamental relations 
\begin{center}
$(T_{0}-u_{1})(T_{0}-u_{2})\cdots (T_{0}-u_{e})=0,$\\
$(T_{i}-q)(T_{i}+1)=0\ \ \ \ \ \ (1\leq i\leq n-1),$\\
$T_{0}T_{1}T_{0}T_{1}=T_{1}T_{0}T_{1}T_{0},$\\
$T_{i}T_{j}=T_{j}T_{i}\ \ \ \ \ \ (|i-j|\geq 2),$\\ 
$T_{i}T_{i+1}T_{i}=T_{i+1}T_{i}T_{i+1}\ \ \ \ \ \ (1\leq i\leq n-2)$.
\end{center}
\end{defi}
\begin{rem} When $e=1$ and $u_{1}=1$, the corresponding Ariki-Koike algebra of type $G(1,1,n)$ is just the Iwahori-Hecke algebra of type $A_{n-1}$. Similarly,
when $e=2$ and $u_{1}=-1$, we have the Iwahori-Hecke algebra of type $B_{n}$ with unequal parameters. 
\end{rem}

For $i\in \{ 0,1,\cdots n-1\}$ we put $L_{i}=T_{i}\cdots T_{1}T_{0}T_{1}^{-1}\cdots T_{i}^{-1}$. We also define an element $T_{w}$ ($w\in \mathfrak{S}_{n}$) as
follows: let $w=s_{i_{1}}\cdots s_{i_{\ell}}$ be a reduced expression of $w$ where $s_{i}$ is the permutation which transposes $i$ and $i+1$. Then we put
$T_{w}=T_{i_{1}}\cdots T_{i_{\ell}}$. By [AK, Theorem 3.10] the set 
\begin{center}
$\{\ L_{0}^{e_{0}}L_{1}^{e_{1}}\cdots L_{n-1}^{e_{n-1}}T_{w}\ |\ 0\leq e_{i}\leq e-1\ (0\leq i\leq n-1),\ w\in \mathfrak{S}_{n}\ \}$
\end{center}
forms a $\Bbbk$-basis of $H_{e,n}$.

Note that for each positive integer $e$ we have the following inductive system:
\begin{center}
$H_{e,1}\subset H_{e,2}\subset \cdots\subset H_{e,n}\subset \cdots$.
\end{center}
We denote by $H_{e}$ the inductive limit of the inductive system, i.e., the union
\begin{center}
$\displaystyle \bigcup_{n=1}^{\infty}H_{e,n}$. 
\end{center}

Next, let us recall the definition of multi-Young tableaux and related notions.
\begin{defi}
An \textbf{$\bm{e}$-Young diagram} of total size $n$ is an $e$-tuple $\bm{\lambda }=(\lambda _{1},\cdots,\lambda _{e})$
consisting of sequences $\lambda _{i}=(\lambda _{i,1},\cdots,\lambda _{i,p(i)})$ of positive integers which satisfies the following 
two conditions:
\begin{itemize}
\item[(1)] $\displaystyle \sum_{\stackrel{\scriptstyle 1\leq i\leq n}{ 1\leq j\leq p(i)}}\lambda _{i,j}=n$,
\item[(2)] $\lambda _{i,1}\geq \lambda _{i,2}\geq \cdots \geq \lambda_{i,p(i)}$ for all $i\in \{ 1,\cdots, e\}$. 
\end{itemize} 
\end{defi}
\begin{rem}
In the definition of the multi-Young diagram, we allow the case that $\lambda _{i}=\emptyset$ for some $i$.
\end{rem}

Obviuously we can regard an $e$-Young diagram as an $e$-tuple of Young diagrams.
\begin{exm}
The $4$-Young diagram $\bm{\lambda }=((4,3,1),\emptyset, (3,1,1),\emptyset)$ of total size $13$ is 
identified with the following 4-tuple of Young diagrams:
\begin {figure} [h]
\[ \left( \
\vcenter {\xy 0;/r.1pc/:
(0,0);(40,0)**@{-},
(0,-10);(40,-10)**@{-},
(0,-20);(10,-20)**@{-},
(0,0);(0,-20)**@{-},
(10,0);(10,-20)**@{-},
(20,0);(20,-10)**@{-},
(30,0);(30,-10)**@{-},
(40,0);(40,-10)**@{-},
(20,-10);(20,-20)**@{-},
(30,-10);(30,-20)**@{-},
(10,-20);(20,-20)**@{-},
(20,-20);(30,-20)**@{-},
(0,-30);(10,-30)**@{-},
(10,-20);(10,-30)**@{-},
(0,-20);(0,-30)**@{-}
\endxy}\ , \ \bm{\emptyset} \ ,\
 \vcenter {\xy 0;/r.1pc/:
(0,0);(30,0)**@{-},
(0,-10);(30,-10)**@{-},
(0,-20);(10,-20)**@{-},
(0,0);(0,-20)**@{-},
(10,0);(10,-20)**@{-},
(20,0);(20,-10)**@{-},
(30,0);(30,-10)**@{-},
(10,-20);(10,-30)**@{-},
(0,-20);(0,-30)**@{-},
(0,-30);(10,-30)**@{-},

\endxy}\ , \ \bm{\emptyset} \ \right) \
\]
\end{figure}
\end{exm}

Let $a,b$ be boxes in an $e$-Young diagram $\bm{\lambda}$. Then the \textbf{content} $c(a;\bm{\lambda })$ of $a$ is
the difference  
\begin{center}
$(\text{the column number of}\ a)-(\text{the row number of}\ a)$ 
\end{center}
and the \textbf{axial distance} $r(a,b)$ from $a$ to $b$ is the difference $c(b;\bm{\lambda })-c(a;\bm{\lambda })$.

\begin{defi}
Let $\bm{\lambda} $ be an $e$-Young diagram. A \textbf{standard $\bm{e}$-tableau} $\bold{T}$ of shape $\bm{\lambda }$ is a pair of 
an $e$-Young diagram and an ordering on the boxes by $\{1,2,\cdots,n\}$ which satisfies the following condition:
in each Young diagram the written numbers are increasing from left to right and from top to bottom.
\end{defi}
\begin{exm}
Let $\bm{\lambda }$ be the $4$-Young diagram presented in Example 2.5. Then in the following two figures
the left one is a standard $4$-tableau of shape $\bm{\lambda }$ but the right one is not.
\begin {figure} [h]
\[ \left( \
\vcenter {\xy 0;/r.1pc/:
(0,0);(40,0)**@{-},
(0,-10);(40,-10)**@{-},
(0,-20);(10,-20)**@{-},
(0,0);(0,-20)**@{-},
(10,0);(10,-20)**@{-},
(20,0);(20,-10)**@{-},
(30,0);(30,-10)**@{-},
(40,0);(40,-10)**@{-},
(20,-10);(20,-20)**@{-},
(30,-10);(30,-20)**@{-},
(10,-20);(20,-20)**@{-},
(20,-20);(30,-20)**@{-},
(0,-30);(10,-30)**@{-},
(10,-20);(10,-30)**@{-},
(0,-20);(0,-30)**@{-},
(5,-5)*{2},
(15,-5)*{5},
(25,-5)*{6},
(35,-5)*{10},
(5,-15)*{8},
(5,-25)*{9},
(15,-15)*{11},
(25,-15)*{12},
\endxy}\ , \ \bm{\emptyset} \ ,\
 \vcenter {\xy 0;/r.1pc/:
(0,0);(30,0)**@{-},
(0,-10);(30,-10)**@{-},
(0,-20);(10,-20)**@{-},
(0,0);(0,-20)**@{-},
(10,0);(10,-20)**@{-},
(20,0);(20,-10)**@{-},
(30,0);(30,-10)**@{-},
(10,-20);(10,-30)**@{-},
(0,-20);(0,-30)**@{-},
(0,-30);(10,-30)**@{-},
(5,-5)*{1},
(5,-15)*{3},
(5,-25)*{13},
(15,-5)*{4},
(25,-5)*{7},
\endxy}\ , \ \bm{\emptyset} \ \right) ,\
 \left( \
\vcenter {\xy 0;/r.1pc/:
(0,0);(40,0)**@{-},
(0,-10);(40,-10)**@{-},
(0,-20);(10,-20)**@{-},
(0,0);(0,-20)**@{-},
(10,0);(10,-20)**@{-},
(20,0);(20,-10)**@{-},
(30,0);(30,-10)**@{-},
(40,0);(40,-10)**@{-},
(20,-10);(20,-20)**@{-},
(30,-10);(30,-20)**@{-},
(10,-20);(20,-20)**@{-},
(20,-20);(30,-20)**@{-},
(0,-30);(10,-30)**@{-},
(10,-20);(10,-30)**@{-},
(0,-20);(0,-30)**@{-},
(5,-5)*{2},
(15,-5)*{5},
(25,-5)*{6},
(35,-5)*{10},
(5,-15)*{11},
(5,-25)*{9},
(15,-15)*{8},
(25,-15)*{12},
\endxy}\ , \ \bm{\emptyset} \ ,\
 \vcenter {\xy 0;/r.1pc/:
(0,0);(30,0)**@{-},
(0,-10);(30,-10)**@{-},
(0,-20);(10,-20)**@{-},
(0,0);(0,-20)**@{-},
(10,0);(10,-20)**@{-},
(20,0);(20,-10)**@{-},
(30,0);(30,-10)**@{-},
(10,-20);(10,-30)**@{-},
(0,-20);(0,-30)**@{-},
(0,-30);(10,-30)**@{-},
(5,-5)*{1},
(5,-15)*{3},
(5,-25)*{13},
(15,-5)*{4},
(25,-5)*{7},
\endxy}\ , \ \bm{\emptyset} \ \right) \
\]
\end{figure}
\end{exm}
For an $e$-Young diagram $\bm{\lambda}$, we denote by $\text{Std}(\bm{\lambda })$ the set of standard $e$-tableaux of shape $\bm{\lambda }$.
For an integer $k$ and an indeterminate $y$ we define a Laurent polynomial $\Delta (k,y)$ and a matrix $M(k,y)$ as follows:
\begin{center}
$\Delta(k,y)=1-q^{k}y$,\ \ \ \ 
$\displaystyle M(k,y)=\frac{1}{\Delta(k,y)}\begin{bmatrix}
q-1 &  \Delta(k+1,y)  \\
q\Delta(k-1,y) & -q^{k}y(q-1)
\end{bmatrix}.$
\end{center} 
Finally, for a standard $e$-tableau $\bold{T}$ we define the number $\tau (i)$ so that $i$ is written in the $\tau(i)$-th Young diagram of $\bold{T}$.
For an $e$-Young diagram $\bm{\lambda }$ of total size $n$ we denote by $V(\bm{\lambda })$ the finite dimensional vector space
\begin{center}
$\displaystyle V(\bm{\lambda })=\bigoplus_{\bold{T}\in \text{Std}(\bm{\lambda })}\Bbbk \bold{t}$
\end{center}
where $\bold{t}$ is the symbol corresponding to a standard $e$-tableau $\bold{T}$.

Following [AK] we define a representation of $H_{e,n}$ on $V(\bm{\lambda})$ as follows:
\begin{itemize}
\item[(1)] $T_{0}\bold{t}=u_{\tau (1)}\bold{t}$,
\item[(2)] For $i\in \{ 1,\cdots,n-1\}$ we define $T_{i}\bold{t}$ as follows:\\
(2-1) We define $T_{i}\bold{t}=q\bold{t}$, if $i$ and $i+1$ are placed as
\[ \vcenter {\xy 0;/r.1pc/:
(0,0);(40,0)**@{-},
(0,0);(0,-20)**@{-},
(0,-20);(40,-20)**@{-},
(40,0);(40,-20)**@{-},
(20,0);(20,-20)**@{-},
(10,-10)*{i},
(30,-10)*{i+1}
\endxy}\]
(2-2) We define $T_{i}\bold{t}=-\bold{t}$, if $i$ and $i+1$ are placed as 
\[ \vcenter {\xy 0;/r.1pc/:
(0,0);(20,0)**@{-},
(0,0);(0,-40)**@{-},
(0,-20);(20,-20)**@{-},
(20,0);(20,-40)**@{-},
(0,-40);(20,-40)**@{-},
(10,-10)*{i},
(10,-30)*{i+1}
\endxy}\]
(2-3) In the other case we define $T_{i}\bold{t}$ by the following:
\begin{center}
$\displaystyle T_{i}\left<\bold{t},\bold{t}'\right>=\left<\bold{t},\bold{t}'\right>M\Bigg(r(i+1,i),\frac{u_{\tau (i)}}{u_{\tau (i+1)}}\Bigg)$.
\end{center} 
Here $\bold{t}'$ is the symbol corresponding to the standard $e$-Young tableau $\bold{T}'$ obtained by transposing $i$ and $i+1$
in $\bold{T}$.
\end{itemize}

\begin{rem}
To unify the above definitions (1),(2-1),(2-2) and (2-3) we introduce the following notation. For the permutation $s_{i}=(i,i+1)$ and
a standard $e$-tableau $\bold{T}$ we define a vector $s_{i}\bold{t}$ in $V(\bm{\lambda})$ by
\begin{eqnarray}
s_{i}\bold{t}=\left\{ \begin{array}{ll}
\bold{t}' & (\text{if}\ \bold{T}'\ \text{is standard}) \\
0 & (\text{else}) \\
\end{array} \right. \nonumber.
\end{eqnarray}
Here, $\bold{T}'$ is the $e$-Young tableau obtained by transposing $i$ and $i+1$ in $\bold{T}$.

Under this notation we have
\[
T_{i}\bold{t}=\frac{(q-1)u_{\tau(i)}}{u_{\tau(i)}-q^{r(i+1,i)}u_{\tau(i+1)}}\bold{t}+
              \frac{q(u_{\tau(i)}-q^{r(i+1,i)-1}u_{\tau(i+1)})}{u_{\tau(i+1)}-q^{r(i+1,i)}u_{\tau(i)}}(s_{i}\bold{t}).
\]
for all $i\in \{0,\ldots,n-1\}$ and a standard $e$-tableau $\bold{T}$.
\end{rem}

The set
\begin{center}
$\{ \ V(\bm{\lambda })\ |\ \bm{\lambda}:e\text{-Young diagram of total size}\ n\ \}$ 
\end{center}
gives a complete list of finite dimensional irreducible representations over $\Bbbk$ up to equivalence.

For an $e$-Young diagram of total size $n$ we denote by $\chi _{\bm{\lambda }}$ the irreducible character corresponding to $\bm{\lambda }$. 

\section{Fourier expansion of Markov traces}
  
In this section we recall a result of Lambropoulou [L2] and of Geck-Iancu-Malle [GIM] concerning the Markov traces of the Ariki-Koike algebras. 
  
Let $z,y_{1},\cdots,y_{e-1}$ be elements of $\Bbbk$. A $\Bbbk$-linear map $\tau :H_{e}\longrightarrow \Bbbk$ is called the \textbf{Markov trace} associated to
$z,y_{1},\cdots,y_{e-1}$ if it satisfies the following properties:
\begin{itemize}
\item[(1)] $\tau (1)=1,$
\item[(2)] $\tau (hh')=\tau (h'h)\ \ (h,h'\in H_{e,n}),$
\item[(3)] $\tau (hT_{i})=z\tau(h)\ \ (h\in H_{e,n}, 1\leq i\leq n),$
\item[(4)] $\tau (hL_{i}^{j})=y_{j}\tau(h)\ \ (h\in H_{e,n}, 0\leq i\leq n-1,\ 1\leq j\leq e-1 ).$
\end{itemize}
  
By [L2] for fixed $e,z,y_{1},\cdots,y_{e-1}$ the Markov trace of $H_{e}$ exists uniquely. 

Since the Markov trace $\tau $ satisfies (2), the restriction $\tau |_{H_{e,n}}$ is written as the $\Bbbk$-linear combination of irreducible characters. 
Geck-Iancu-Malle [GIM] determined the coefficients in this expression.

To state their result we prepare some notations. Let $\bm{\lambda }=(\lambda _{1},\cdots,\lambda _{e})$ be an $e$-Young diagram
of total size $n$. By adding some zeros we regard $\lambda _{1}$ and $\lambda _{p}\ (p\in \{ 2,\cdots,e\})$ as the sequences
\begin{center}
$\lambda _{1}=(\lambda_{1,1},\cdots,\lambda _{1,n},\lambda _{1,n+1}),$\\
$\lambda _{p}=(\lambda_{p,1},\cdots,\lambda _{p,n})$
\end{center}
of length $n+1$ and $n$ respectively. We define finite sets $\bold{A}_{1},\cdots,\bold{A}_{e}$ as follows:
\begin{center}
$\bold{A}_{1}=\{ \alpha _{1,i}:=\lambda _{1,i}+n-i+1\ |\ 1\leq i\leq n+1 \ \},$\\
$\bold{A}_{p}=\{ \alpha_{p,i}:=\lambda _{p,i}+n-i\ |\ 1\leq i\leq n \ \}\ \ \ \ \ (2\leq p\leq n).$
\end{center}
Finally, for a non-negative integer $d$ we denote by $\sigma _{d}$ the $d$-th fundamental symmetric polynomial with respect to $u_{1},\cdots,u_{e}$.
Here, we understand $\sigma _{0}=1$. Under these notations, we define $\widehat{D}^{\bm{\lambda }}(q,u_{1},\cdots,u_{e})$ and $R^{\bm{\lambda}}(z,y_{1},\cdots,y_{e-1})$ as follows:\\
\linebreak
$\widehat{D}^{\bm{\lambda }}(q,u_{1},\cdots,u_{e})$\\
$\displaystyle 
=\frac{\displaystyle (-1)^{\binom{e}{2}\binom{n}{2}+n(e-1)}\prod_{1\leq k\leq l\leq e}
\prod_{\stackrel{\scriptstyle (\alpha ,\alpha ')\in \bold{A}_{k}\times \bold{A}_{l}}{\alpha>\alpha' \ (\text{if}\ k=l)}}(q^{\alpha }u_{k}-q^{\alpha'}u_{l})\prod_{k=1}^{e}u_{k}^{n}}                                
{\displaystyle  q^{f(n,e)}\prod_{k,l=1}^{e}\prod_{\alpha \in \bold{A}_{k}}\prod_{h=1}^{\alpha }(q^{h}u_{k}-u_{l})\prod_{1\leq k<l\leq e}(u_{k}-u_{l}) }$,\\
\linebreak
$\displaystyle R^{\bm{\lambda}}(z,y_{1},\cdots,y_{e-1})=\prod_{k=1}^{e}\prod_{x\in \lambda _{k}}\Bigg( z(1-q^{c(x)})
\prod_{\stackrel{\scriptstyle 1\leq l\leq e}{l\not= k}}(q^{c(x)}u_{k}-u_{l})$\\
$\hspace{20mm} \displaystyle +(1-q)( \prod_{i=1}^{e-1}(q^{(e-i)c(x)}u_{k}^{e-i-1}\sum_{j=1}^{i}(-1)^{i-j}y_{j}\sigma _{i-j})+(-1)^{e-1}
\prod_{\stackrel{\scriptstyle 1\leq l\leq e}{l\not= k}}u_{l}) \Bigg)$.\\
\linebreak
Here, $f(n,e)$ is given by
\begin{center}
$\displaystyle \sum_{i=1}^{n-1}\binom{ei+1}{2}=\binom{e(n-1)+1}{2}+\binom{e(n-2)+1}{2}+\cdots+\binom{e+1}{2}$.
\end{center}
Finally, we put $D^{\bm{\lambda}}(q,u_{1},u_{2},\cdots,u_{e})=\widehat{D}^{\bm{\lambda}}(q,q^{-1}u_{1},u_{2},\cdots,u_{e})$.
\begin{lem}{ ([GIM])}
We have
\begin{center}
$\displaystyle \tau |_{H_{e,n}}=\sum_{\bm{\lambda}}(-1)^{en}\Bigg(\prod_{k=1}^{e}u_{k}^{|\lambda _{k}|-n}\Bigg)
D^{\bm{\lambda}}(q,u_{1},u_{2},\cdots,u_{e})R^{\bm{\lambda}}(z,y_{1},\cdots,y_{e-1})\chi _{\bm{\lambda}}$
\end{center}
where $\bm{\lambda}$ runs through all $e$-Young diagrams of total size $n$ and $|\lambda _{k}|$ stands for the number of boxes in $\lambda _{k}$.
\end{lem}
\begin{rem}
We see easily by induction on $n$ that
\begin{center}
$\displaystyle f(n,e)=\frac{1}{12}en(n-1)(2en-e+3)$.
\end{center}
\end{rem}

In the rest of this paper, we set
\begin{center}
$\displaystyle C(e,k,x)=\prod_{i=1}^{e-1}(q^{(e-i)c(x)}u_{k}^{e-i-1}\sum_{j=1}^{i}(-1)^{i-j}y_{j}\sigma _{i-j})+(-1)^{e-1}\prod_{\stackrel{\scriptstyle 1\leq l\leq e}{l\not= k}}u_{l}$.
\end{center}

\section{Mixed link theory}
In this section we collect some fundamental notions and facts from Lambropoulou's mixed link theory.
For the details we refer [L1] and [L2]. 

\begin{defi}
\begin{itemize}
\item[(1)] A \textbf{mixed link} is an embedding of a disjoint union of finitely many circles into the solid torus.
\item[(2)] Two mixed links are said to be \textbf{equivalent} if these are joined each other by an ambient isotropy of the solid torus.
\end{itemize}
\end{defi}
By considering the canonical genus 1 Heegaard decomposition, we regard the solid torus as the complement of the interior of a solid torus: 
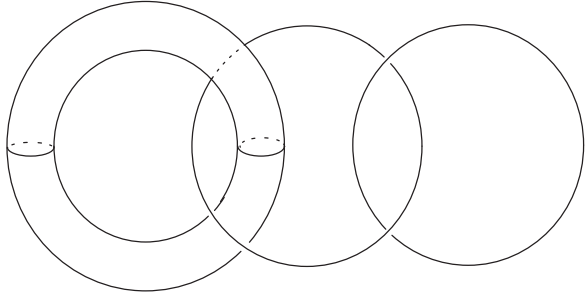
\begin{center}
\unitlength 0.1in
\begin{picture}( 30.1600, 15.1400)(  4.0400,-19.2000)
%
{\color[named]{Black}{%
\special{pn 8}%
\special{ar 1130 1164 726 758  0.8137996 6.2831853}%
\special{ar 1130 1164 726 758  0.0000000 0.7259561}%
}}%
%
{\color[named]{Black}{%
\special{pn 8}%
\special{ar 1130 1164 480 502  0.8188646 6.2831853}%
\special{ar 1130 1164 480 502  0.0000000 0.6168629}%
}}%
%
{\color[named]{Black}{%
\special{pn 8}%
\special{ar 528 1170 122 46  6.1517166 6.2831853}%
\special{ar 528 1170 122 46  0.0000000 3.1415927}%
}}%
%
{\color[named]{Black}{%
\special{pn 8}%
\special{ar 1734 1164 122 52  6.2831853 6.2831853}%
\special{ar 1734 1164 122 52  0.0000000 3.1415927}%
}}%
%
{\color[named]{Black}{%
\special{pn 8}%
\special{ar 528 1176 122 32  3.1415927 3.2994874}%
\special{ar 528 1176 122 32  3.7731716 3.9310663}%
\special{ar 528 1176 122 32  4.4047505 4.5626453}%
\special{ar 528 1176 122 32  5.0363295 5.1942242}%
\special{ar 528 1176 122 32  5.6679084 5.8258032}%
}}%
%
{\color[named]{Black}{%
\special{pn 8}%
\special{ar 1740 1164 116 44  3.0152845 3.1662279}%
\special{ar 1740 1164 116 44  3.6190581 3.7700015}%
\special{ar 1740 1164 116 44  4.2228317 4.3737751}%
\special{ar 1740 1164 116 44  4.8266053 4.9775487}%
\special{ar 1740 1164 116 44  5.4303789 5.5813223}%
\special{ar 1740 1164 116 44  6.0341525 6.1850959}%
}}%
%
{\color[named]{Black}{%
\special{pn 8}%
\special{ar 1974 1164 604 630  4.1375807 5.2887642}%
}}%
%
{\color[named]{Black}{%
\special{pn 8}%
\special{ar 1974 1164 602 628  3.7402628 3.7597908}%
\special{ar 1974 1164 602 628  3.8183750 3.8379031}%
\special{ar 1974 1164 602 628  3.8964873 3.9160154}%
\special{ar 1974 1164 602 628  3.9745996 3.9941277}%
\special{ar 1974 1164 602 628  4.0527119 4.0722400}%
}}%
%
{\color[named]{Black}{%
\special{pn 8}%
\special{ar 1974 1164 602 630  5.5274344 6.2831853}%
\special{ar 1974 1164 602 630  0.0000000 3.7402628}%
}}%
%
{\color[named]{Black}{%
\special{pn 8}%
\special{ar 2818 1158 604 630  2.3662613 6.2831853}%
\special{ar 2818 1158 604 630  0.0000000 2.2927503}%
}}%
%
{\color[named]{Black}{%
\special{pn 8}%
\special{ar 1974 1164 602 630  5.2842267 5.4482993}%
}}%
%
{\color[named]{Black}{%
\special{pn 8}%
\special{ar 1130 1158 488 508  0.5700854 0.7406949}%
}}%
\end{picture}%

\end{center}
\begin{defi}
(1) A a mixed braid with $n$-strands is an $n$-tuple $(p_{1},\dots,p_{n})$ of embeddings of the closed interval $[0,1]$ into the closure $\text{cl}(I^{3}\setminus Cyn)$,
where 
\begin{center}
$\displaystyle Cyn:=\bigg\{(x,y,z)\in I^{3} \bigg| \bigg( x-\frac{1}{4}\bigg)^{2}+\bigg( y-\frac{1}{2(n+1)}\bigg)^{2}\leq \bigg(\frac{1}{4(n+1)}\bigg)^{2}\bigg\}$,
\end{center}
such that
\begin{itemize}
\item the curves $p_{i}(I)$ and $p_{j}(I)$ do not intersect if $i\not= j$,
\item there exists a permutation $\sigma $ of $\{ 1,\dots,n\}$ such that 
\begin{center}
$\displaystyle p_{i}(0)=\bigg(\frac{1}{2},\frac{i}{n+1},0\bigg),\ p_{i}(1)=\bigg(\frac{1}{2},\frac{\sigma(i)}{n+1},1\bigg)$
\end{center}
for all $i\in \{ 1,\dots, n\}$,
\item the point $p_{i}(t)$ lies in the interior of $\text{cl}(I^{3}\setminus Cyn)$ if $t\in (0,1)$,
\item the 3rd coordinate of $p_{i}(t)$ is increasing with respect to $t$ for all $i\in \{ 1,\dots, n\}$.
\end{itemize}
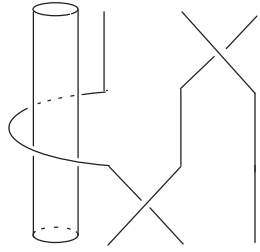
\begin{center}
\unitlength 0.1in
\begin{picture}( 12.9900, 12.7000)( 12.1100,-22.1000)
%
{\color[named]{Black}{%
\special{pn 8}%
\special{ar 1456 974 118 34  0.8934099 6.2831853}%
\special{ar 1456 974 118 34  0.0000000 0.7853982}%
}}%
%
{\color[named]{Black}{%
\special{pn 8}%
\special{ar 1456 2156 118 40  6.1049178 6.2831853}%
\special{ar 1456 2156 118 40  0.0000000 3.5204240}%
}}%
%
{\color[named]{Black}{%
\special{pn 8}%
\special{ar 1456 2156 118 40  2.8250234 2.9778896}%
\special{ar 1456 2156 118 40  3.4364884 3.5893546}%
\special{ar 1456 2156 118 40  4.0479533 4.2008196}%
\special{ar 1456 2156 118 40  4.6594183 4.8122845}%
\special{ar 1456 2156 118 40  5.2708833 5.4237495}%
\special{ar 1456 2156 118 40  5.8823482 6.0352145}%
}}%
%
{\color[named]{Black}{%
\special{pn 8}%
\special{pa 1730 2210}%
\special{pa 2112 1798}%
\special{fp}%
}}%
%
{\color[named]{Black}{%
\special{pn 8}%
\special{pa 2124 2204}%
\special{pa 1950 2014}%
\special{fp}%
}}%
%
{\color[named]{Black}{%
\special{pn 8}%
\special{pa 1736 1798}%
\special{pa 1904 1974}%
\special{fp}%
}}%
%
{\color[named]{Black}{%
\special{pn 8}%
\special{pa 2500 2196}%
\special{pa 2500 1792}%
\special{fp}%
}}%
%
{\color[named]{Black}{%
\special{pn 8}%
\special{ar 1918 1596 706 204  1.8311145 3.6687386}%
}}%
%
{\color[named]{Black}{%
\special{pn 8}%
\special{ar 1918 1596 714 196  4.2429049 4.4259287}%
}}%
%
{\color[named]{Black}{%
\special{pn 8}%
\special{ar 1918 1590 720 190  3.8352724 3.8617041}%
\special{ar 1918 1590 720 190  3.9409992 3.9674309}%
\special{ar 1918 1590 720 190  4.0467261 4.0731578}%
\special{ar 1918 1590 720 190  4.1524530 4.1788847}%
}}%
%
{\color[named]{Black}{%
\special{pn 8}%
\special{pa 1340 988}%
\special{pa 1340 1690}%
\special{fp}%
}}%
%
{\color[named]{Black}{%
\special{pn 8}%
\special{pa 1340 1744}%
\special{pa 1340 2156}%
\special{fp}%
}}%
%
{\color[named]{Black}{%
\special{pn 8}%
\special{pa 1574 982}%
\special{pa 1574 1752}%
\special{fp}%
}}%
%
{\color[named]{Black}{%
\special{pn 8}%
\special{pa 1574 2156}%
\special{pa 1574 1798}%
\special{fp}%
}}%
%
{\color[named]{Black}{%
\special{pn 8}%
\special{pa 1722 1400}%
\special{pa 1722 1400}%
\special{fp}%
}}%
%
{\color[named]{Black}{%
\special{pn 8}%
\special{pa 1710 1400}%
\special{pa 1710 988}%
\special{fp}%
}}%
%
{\color[named]{Black}{%
\special{pn 8}%
\special{pa 2500 1414}%
\special{pa 2500 1826}%
\special{fp}%
}}%
%
{\color[named]{Black}{%
\special{pn 8}%
\special{pa 2500 1414}%
\special{pa 2118 988}%
\special{fp}%
}}%
%
{\color[named]{Black}{%
\special{pn 8}%
\special{pa 2112 1792}%
\special{pa 2112 1394}%
\special{fp}%
}}%
%
{\color[named]{Black}{%
\special{pn 8}%
\special{pa 2110 1390}%
\special{pa 2286 1222}%
\special{fp}%
}}%
%
{\color[named]{Black}{%
\special{pn 8}%
\special{pa 2510 1010}%
\special{pa 2350 1180}%
\special{fp}%
}}%
\end{picture}%

\end{center}
(2) Two mixed braids with $n$-strands are said to be \textbf{equivalent} if they are joined by an ambient isotopy of $\text{cl}(I^{3}\setminus Cyn)$ which
fixes the boundary.
\end{defi}

As in the usual braid theory the set of equivalence classes of  mixed links with $n$-strands is equipped with the natural group structure.

Let $B_{\text{aff},n}$ be the $n$-th affine braid group, namely, the group with generators 
$t_{0},t_{1},\cdots,t_{n-1}$ 
satisfying fundamental relations 
\begin{center}
$t_{0}t_{1}t_{0}t_{1}=t_{1}t_{0}t_{1}t_{0},$\\
$t_{i}t_{j}=t_{j}t_{i}\ \ \ \ \ \ (|i-j|\geq 2),$\\
$t_{i}t_{i+1}t_{i}=t_{i+1}t_{i}t_{i+1}\ \ \ \ \ \  (1\leq i\leq n-2)$.
\end{center}
\begin{lem}
The group of equivelence classes of mixed braids with $n$-strands is isomorphic to $B_{\text{aff},n}$ as groups. Here, $t_{0}$ corresponds to
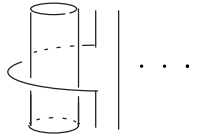
\begin{center}
\unitlength 0.1in
\begin{picture}( 10.7000,  6.8000)( 16.9000,-19.5000)
%
{\color[named]{Black}{%
\special{pn 8}%
\special{pa 2060 1300}%
\special{pa 2060 1710}%
\special{fp}%
}}%
%
{\color[named]{Black}{%
\special{pn 8}%
\special{pa 2060 1790}%
\special{pa 2060 1900}%
\special{fp}%
\special{pa 2060 1900}%
\special{pa 2060 1900}%
\special{fp}%
}}%
%
{\color[named]{Black}{%
\special{pn 8}%
\special{pa 1810 1750}%
\special{pa 1810 1910}%
\special{fp}%
}}%
%
{\color[named]{Black}{%
\special{pn 8}%
\special{ar 1930 1910 130 40  6.0418298 6.2831853}%
\special{ar 1930 1910 130 40  0.0000000 3.1415927}%
}}%
%
{\color[named]{Black}{%
\special{pn 8}%
\special{ar 1930 1910 130 40  3.4021950 3.5433715}%
\special{ar 1930 1910 130 40  3.9669009 4.1080774}%
\special{ar 1930 1910 130 40  4.5316068 4.6727833}%
\special{ar 1930 1910 130 40  5.0963127 5.2374892}%
\special{ar 1930 1910 130 40  5.6610186 5.8021950}%
\special{ar 1930 1910 130 40  6.2257245 6.3669009}%
}}%
%
{\color[named]{Black}{%
\special{pn 8}%
\special{ar 2180 1630 490 100  1.6115900 3.6911527}%
}}%
%
{\color[named]{Black}{%
\special{pn 8}%
\special{ar 2160 1590 490 100  4.5505515 4.6753002}%
}}%
%
{\color[named]{Black}{%
\special{pn 8}%
\special{pa 2150 1500}%
\special{pa 2150 1310}%
\special{fp}%
}}%
%
{\color[named]{Black}{%
\special{pn 8}%
\special{pa 2150 1730}%
\special{pa 2150 1920}%
\special{fp}%
}}%
%
{\color[named]{Black}{%
\special{pn 8}%
\special{ar 2150 1600 430 110  3.8654249 3.9098693}%
\special{ar 2150 1600 430 110  4.0432027 4.0876471}%
\special{ar 2150 1600 430 110  4.2209804 4.2654249}%
\special{ar 2150 1600 430 110  4.3987582 4.4432027}%
}}%
%
{\color[named]{Black}{%
\special{pn 8}%
\special{pa 1810 1660}%
\special{pa 1810 1320}%
\special{fp}%
}}%
%
{\color[named]{Black}{%
\special{pn 8}%
\special{pa 2270 1310}%
\special{pa 2270 1930}%
\special{fp}%
}}%
%
{\color[named]{Black}{%
\special{pn 4}%
\special{sh 1}%
\special{ar 2390 1600 6 6 0  6.28318530717959E+0000}%
\special{sh 1}%
\special{ar 2510 1600 6 6 0  6.28318530717959E+0000}%
\special{sh 1}%
\special{ar 2630 1600 6 6 0  6.28318530717959E+0000}%
\special{sh 1}%
\special{ar 2630 1600 6 6 0  6.28318530717959E+0000}%
}}%
%
{\color[named]{Black}{%
\special{pn 8}%
\special{pa 2760 1320}%
\special{pa 2760 1920}%
\special{fp}%
}}%
%
{\color[named]{Black}{%
\special{pn 8}%
\special{ar 1930 1300 120 30  1.0584069 6.2831853}%
\special{ar 1930 1300 120 30  0.0000000 0.7853982}%
}}%
%
{\color[named]{Black}{%
\special{pn 8}%
\special{pa 2060 1890}%
\special{pa 2060 1760}%
\special{fp}%
}}%
%
{\color[named]{Black}{%
\special{pn 8}%
\special{pa 1810 1910}%
\special{pa 1810 1730}%
\special{fp}%
}}%
\end{picture}%

\end{center}
and $t_{i}\ (1\leq i\leq n-1)$ corresponds to
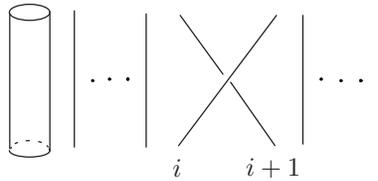
\begin{center}
\unitlength 0.1in
\begin{picture}( 19.5100,  7.9700)(  7.1900,-13.8800)
%
{\color[named]{Black}{%
\special{pn 8}%
\special{ar 824 632 106 42  0.0000000 6.2831853}%
}}%
%
{\color[named]{Black}{%
\special{pn 8}%
\special{pa 720 640}%
\special{pa 720 1346}%
\special{fp}%
}}%
%
{\color[named]{Black}{%
\special{pn 8}%
\special{pa 930 640}%
\special{pa 930 1356}%
\special{fp}%
}}%
%
{\color[named]{Black}{%
\special{pn 8}%
\special{pa 1058 1338}%
\special{pa 1058 624}%
\special{fp}%
}}%
%
{\color[named]{Black}{%
\special{pn 8}%
\special{ar 824 1356 106 34  6.0137127 6.2831853}%
\special{ar 824 1356 106 34  0.0000000 3.1415927}%
}}%
%
{\color[named]{Black}{%
\special{pn 8}%
\special{ar 824 1346 98 42  3.1415927 3.3142545}%
\special{ar 824 1346 98 42  3.8322401 4.0049020}%
\special{ar 824 1346 98 42  4.5228876 4.6955495}%
\special{ar 824 1346 98 42  5.2135351 5.3861970}%
\special{ar 824 1346 98 42  5.9041826 6.0768445}%
}}%
%
{\color[named]{Black}{%
\special{pn 8}%
\special{pa 1434 640}%
\special{pa 1434 1338}%
\special{fp}%
}}%
%
{\color[named]{Black}{%
\special{pn 4}%
\special{sh 1}%
\special{ar 1242 982 6 6 0  6.28318530717959E+0000}%
\special{sh 1}%
\special{ar 1154 990 6 6 0  6.28318530717959E+0000}%
\special{sh 1}%
\special{ar 1338 982 6 6 0  6.28318530717959E+0000}%
\special{sh 1}%
\special{ar 1338 982 6 6 0  6.28318530717959E+0000}%
}}%
%
{\color[named]{Black}{%
\special{pn 8}%
\special{pa 1604 1330}%
\special{pa 2116 648}%
\special{fp}%
}}%
%
{\color[named]{Black}{%
\special{pn 8}%
\special{pa 2108 1330}%
\special{pa 1876 998}%
\special{fp}%
}}%
%
{\color[named]{Black}{%
\special{pn 8}%
\special{pa 1612 648}%
\special{pa 1836 956}%
\special{fp}%
}}%
%
{\color[named]{Black}{%
\special{pn 8}%
\special{pa 2254 1322}%
\special{pa 2254 648}%
\special{fp}%
}}%
%
{\color[named]{Black}{%
\special{pn 4}%
\special{sh 1}%
\special{ar 2350 982 6 6 0  6.28318530717959E+0000}%
\special{sh 1}%
\special{ar 2454 990 6 6 0  6.28318530717959E+0000}%
\special{sh 1}%
\special{ar 2558 990 6 6 0  6.28318530717959E+0000}%
\special{sh 1}%
\special{ar 2558 990 6 6 0  6.28318530717959E+0000}%
}}%
%
{\color[named]{Black}{%
\special{pn 8}%
\special{pa 2670 666}%
\special{pa 2670 1322}%
\special{fp}%
}}%
\put(15.9500,-14.4600){\makebox(0,0){$i$}}%
\put(21.0000,-14.4600){\makebox(0,0){$i+1$}}%
\end{picture}%

\end{center}

\end{lem}

As in the usual braid theory, we can consider the closure of a mixed braid. 
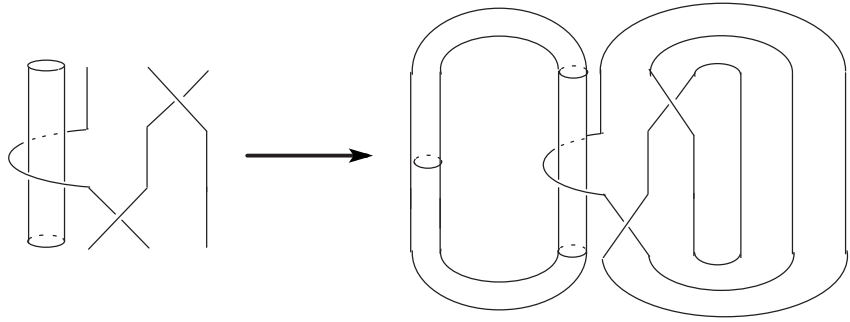
\begin{center}
\unitlength 0.1in
\begin{picture}( 43.8900, 16.4100)( 10.6100,-21.7000)
%
{\color[named]{Black}{%
\special{pn 8}%
\special{ar 1258 856 96 28  0.8982175 6.2831853}%
\special{ar 1258 856 96 28  0.0000000 0.7853982}%
}}%
%
{\color[named]{Black}{%
\special{pn 8}%
\special{ar 1258 1776 96 32  6.0857897 6.2831853}%
\special{ar 1258 1776 96 32  0.0000000 3.5654603}%
}}%
%
{\color[named]{Black}{%
\special{pn 8}%
\special{ar 1258 1776 96 32  2.8198421 3.0103183}%
\special{ar 1258 1776 96 32  3.5817469 3.7722231}%
\special{ar 1258 1776 96 32  4.3436516 4.5341278}%
\special{ar 1258 1776 96 32  5.1055564 5.2960326}%
\special{ar 1258 1776 96 32  5.8674611 6.0579373}%
}}%
%
{\color[named]{Black}{%
\special{pn 8}%
\special{pa 1480 1818}%
\special{pa 1786 1498}%
\special{fp}%
}}%
%
{\color[named]{Black}{%
\special{pn 8}%
\special{pa 1796 1812}%
\special{pa 1654 1666}%
\special{fp}%
}}%
%
{\color[named]{Black}{%
\special{pn 8}%
\special{pa 1482 1498}%
\special{pa 1620 1634}%
\special{fp}%
}}%
%
{\color[named]{Black}{%
\special{pn 8}%
\special{pa 2098 1808}%
\special{pa 2098 1492}%
\special{fp}%
}}%
%
{\color[named]{Black}{%
\special{pn 8}%
\special{ar 1630 1340 568 158  1.8281200 3.6642770}%
}}%
%
{\color[named]{Black}{%
\special{pn 8}%
\special{ar 1630 1340 574 154  4.2414363 4.4218062}%
}}%
%
{\color[named]{Black}{%
\special{pn 8}%
\special{ar 1630 1336 578 146  3.8387358 3.8718850}%
\special{ar 1630 1336 578 146  3.9713325 4.0044817}%
\special{ar 1630 1336 578 146  4.1039292 4.1370784}%
}}%
%
{\color[named]{Black}{%
\special{pn 8}%
\special{pa 1166 866}%
\special{pa 1166 1414}%
\special{fp}%
}}%
%
{\color[named]{Black}{%
\special{pn 8}%
\special{pa 1166 1456}%
\special{pa 1166 1776}%
\special{fp}%
}}%
%
{\color[named]{Black}{%
\special{pn 8}%
\special{pa 1354 862}%
\special{pa 1354 1462}%
\special{fp}%
}}%
%
{\color[named]{Black}{%
\special{pn 8}%
\special{pa 1354 1776}%
\special{pa 1354 1498}%
\special{fp}%
}}%
%
{\color[named]{Black}{%
\special{pn 8}%
\special{pa 1472 1188}%
\special{pa 1472 1188}%
\special{fp}%
}}%
%
{\color[named]{Black}{%
\special{pn 8}%
\special{pa 1472 1188}%
\special{pa 1472 866}%
\special{fp}%
}}%
%
{\color[named]{Black}{%
\special{pn 8}%
\special{pa 2098 1200}%
\special{pa 2098 1518}%
\special{fp}%
}}%
%
{\color[named]{Black}{%
\special{pn 8}%
\special{pa 2098 1200}%
\special{pa 1792 866}%
\special{fp}%
}}%
%
{\color[named]{Black}{%
\special{pn 8}%
\special{pa 1786 1492}%
\special{pa 1786 1182}%
\special{fp}%
}}%
%
{\color[named]{Black}{%
\special{pn 8}%
\special{pa 1786 1178}%
\special{pa 1928 1046}%
\special{fp}%
}}%
%
{\color[named]{Black}{%
\special{pn 8}%
\special{pa 2106 884}%
\special{pa 1978 1018}%
\special{fp}%
}}%
%
{\color[named]{Black}{%
\special{pn 20}%
\special{pa 2310 1328}%
\special{pa 2934 1328}%
\special{fp}%
\special{sh 1}%
\special{pa 2934 1328}%
\special{pa 2866 1308}%
\special{pa 2880 1328}%
\special{pa 2866 1348}%
\special{pa 2934 1328}%
\special{fp}%
\special{pa 2934 1328}%
\special{pa 2934 1328}%
\special{fp}%
}}%
%
{\color[named]{Black}{%
\special{pn 8}%
\special{ar 4012 1828 74 32  6.0553283 6.2831853}%
\special{ar 4012 1828 74 32  0.0000000 3.5464844}%
}}%
%
{\color[named]{Black}{%
\special{pn 8}%
\special{ar 4012 1828 74 32  2.8323937 3.0609652}%
\special{ar 4012 1828 74 32  3.7466795 3.9752509}%
\special{ar 4012 1828 74 32  4.6609652 4.8895366}%
\special{ar 4012 1828 74 32  5.5752509 5.8038223}%
}}%
%
{\color[named]{Black}{%
\special{pn 8}%
\special{pa 4174 1862}%
\special{pa 4408 1530}%
\special{fp}%
}}%
%
{\color[named]{Black}{%
\special{pn 8}%
\special{pa 4416 1856}%
\special{pa 4308 1704}%
\special{fp}%
}}%
%
{\color[named]{Black}{%
\special{pn 8}%
\special{pa 4176 1530}%
\special{pa 4280 1672}%
\special{fp}%
}}%
%
{\color[named]{Black}{%
\special{pn 8}%
\special{pa 4648 1852}%
\special{pa 4648 1526}%
\special{fp}%
}}%
%
{\color[named]{Black}{%
\special{pn 8}%
\special{ar 4296 1376 436 164  1.8295672 3.6627683}%
}}%
%
{\color[named]{Black}{%
\special{pn 8}%
\special{ar 4288 1366 440 160  4.2423788 4.4216187}%
}}%
%
{\color[named]{Black}{%
\special{pn 8}%
\special{ar 4296 1372 444 152  3.8351861 3.8754546}%
\special{ar 4296 1372 444 152  3.9962599 4.0365284}%
\special{ar 4296 1372 444 152  4.1573338 4.1885128}%
}}%
%
{\color[named]{Black}{%
\special{pn 8}%
\special{pa 3940 886}%
\special{pa 3940 1452}%
\special{fp}%
}}%
%
{\color[named]{Black}{%
\special{pn 8}%
\special{pa 3940 1496}%
\special{pa 3940 1828}%
\special{fp}%
}}%
%
{\color[named]{Black}{%
\special{pn 8}%
\special{pa 4086 882}%
\special{pa 4086 1502}%
\special{fp}%
}}%
%
{\color[named]{Black}{%
\special{pn 8}%
\special{pa 4086 1828}%
\special{pa 4086 1538}%
\special{fp}%
}}%
%
{\color[named]{Black}{%
\special{pn 8}%
\special{pa 4168 1210}%
\special{pa 4168 1210}%
\special{fp}%
}}%
%
{\color[named]{Black}{%
\special{pn 8}%
\special{pa 4160 1210}%
\special{pa 4160 878}%
\special{fp}%
}}%
%
{\color[named]{Black}{%
\special{pn 8}%
\special{pa 4648 1222}%
\special{pa 4648 1552}%
\special{fp}%
}}%
%
{\color[named]{Black}{%
\special{pn 8}%
\special{pa 4648 1222}%
\special{pa 4414 876}%
\special{fp}%
}}%
%
{\color[named]{Black}{%
\special{pn 8}%
\special{pa 4408 1526}%
\special{pa 4408 1204}%
\special{fp}%
}}%
%
{\color[named]{Black}{%
\special{pn 8}%
\special{pa 4408 1200}%
\special{pa 4516 1062}%
\special{fp}%
}}%
%
{\color[named]{Black}{%
\special{pn 8}%
\special{pa 4654 896}%
\special{pa 4554 1032}%
\special{fp}%
}}%
%
{\color[named]{Black}{%
\special{pn 8}%
\special{pa 4894 894}%
\special{pa 4894 1862}%
\special{fp}%
}}%
%
{\color[named]{Black}{%
\special{pn 8}%
\special{pa 5164 894}%
\special{pa 5164 1842}%
\special{fp}%
\special{pa 5442 894}%
\special{pa 5442 1842}%
\special{fp}%
}}%
%
{\color[named]{Black}{%
\special{pn 8}%
\special{ar 3256 1360 72 34  6.0942881 6.2831853}%
\special{ar 3256 1360 72 34  0.0000000 3.5329978}%
}}%
%
{\color[named]{Black}{%
\special{pn 8}%
\special{ar 3256 1360 72 34  2.8198421 3.0484135}%
\special{ar 3256 1360 72 34  3.7341278 3.9626992}%
\special{ar 3256 1360 72 34  4.6484135 4.8769850}%
\special{ar 3256 1360 72 34  5.5626992 5.7912707}%
}}%
%
{\color[named]{Black}{%
\special{pn 8}%
\special{pa 3166 892}%
\special{pa 3166 1606}%
\special{fp}%
}}%
%
{\color[named]{Black}{%
\special{pn 8}%
\special{pa 3166 1500}%
\special{pa 3166 1832}%
\special{fp}%
}}%
%
{\color[named]{Black}{%
\special{pn 8}%
\special{pa 3320 890}%
\special{pa 3320 1708}%
\special{fp}%
}}%
%
{\color[named]{Black}{%
\special{pn 8}%
\special{pa 3320 1840}%
\special{pa 3320 1552}%
\special{fp}%
}}%
%
{\color[named]{Black}{%
\special{pn 8}%
\special{ar 3630 884 456 328  3.0867028 6.2831853}%
}}%
%
{\color[named]{Black}{%
\special{pn 8}%
\special{ar 3630 874 308 172  2.9864965 6.2831853}%
}}%
%
{\color[named]{Black}{%
\special{pn 8}%
\special{ar 4018 890 74 32  6.0584741 6.2831853}%
\special{ar 4018 890 74 32  0.0000000 3.5464844}%
}}%
%
{\color[named]{Black}{%
\special{pn 8}%
\special{ar 4018 890 74 32  2.8030370 3.0316084}%
\special{ar 4018 890 74 32  3.7173227 3.9458941}%
\special{ar 4018 890 74 32  4.6316084 4.8601798}%
\special{ar 4018 890 74 32  5.5458941 5.7744655}%
}}%
%
{\color[named]{Black}{%
\special{pn 8}%
\special{ar 3630 1836 308 154  6.2246204 6.2831853}%
\special{ar 3630 1836 308 154  0.0000000 3.2583234}%
}}%
%
{\color[named]{Black}{%
\special{pn 8}%
\special{ar 3630 1836 456 300  6.2831853 6.2831853}%
\special{ar 3630 1836 456 300  0.0000000 3.1723522}%
}}%
%
{\color[named]{Black}{%
\special{pn 8}%
\special{ar 4800 884 642 354  3.0901560 6.2831853}%
}}%
%
{\color[named]{Black}{%
\special{pn 8}%
\special{ar 4792 892 372 192  3.0424563 6.2831853}%
}}%
%
{\color[named]{Black}{%
\special{pn 8}%
\special{ar 4772 902 120 56  2.9617392 6.2831853}%
}}%
%
{\color[named]{Black}{%
\special{pn 8}%
\special{ar 4772 1836 120 64  0.1407302 3.1415927}%
}}%
%
{\color[named]{Black}{%
\special{pn 8}%
\special{ar 4792 1844 380 190  6.2347058 6.2831853}%
\special{ar 4792 1844 380 190  0.0000000 3.0943837}%
}}%
%
{\color[named]{Black}{%
\special{pn 8}%
\special{ar 4810 1854 642 318  6.1975915 6.2831853}%
\special{ar 4810 1854 642 318  0.0000000 3.0854893}%
}}%
\end{picture}%

\end{center}
\begin{lem}
Every mixed link is equivalent to the closure of a mixed braid.
\end{lem}

For a mixed braid $\alpha $ we denote by $\widehat{\alpha }$ the closure of $\alpha $. 
The analogue of the Markov Moves and Alexander's theorem are given by the following lemma.
\begin{lem}
The closures of two mixed braids are equivalent as mixed links if and only if the corresponding mixed braids are joined by a sequence of the following 
two transformations:
\begin{itemize}
\item[(1)] $\alpha \longleftrightarrow \beta \alpha \beta ^{-1}\ \ (\alpha ,\beta \in B_{\text{aff},n})$
\item[(2)] $\alpha \longleftrightarrow \alpha t_{n}^{\pm 1}\ \ (\alpha \in B_{\text{aff},n})$.
\end{itemize}
\end{lem}  

Let us recall the construction of the analogue of the HOMFLYPT polynomial in mixed link theory.
We first introduce new variable $t$ as 
\begin{center}
$\displaystyle t=\frac{1-q+z}{qz}$.
\end{center}
Let $\pi_{n}:B_{\text{aff},n}\longrightarrow H_{e,n}^{\times }$ be the group homomorphism defined by $\pi(t_{i})=T_{i}\ (0\leq i\leq n-1)$.

\begin{defi}
Let $\widehat{\alpha} $ be a mixed link obtained as the closure of a mixed braid $\alpha $ with $n$-strands. Then we define the \textbf{HOMFLYPT polynomial of type $\bm{G(e,1)}$} of $\widehat{\alpha}$
by
\begin{center}
$\displaystyle X_{\widehat{\alpha }}(q,t)=\Bigg[ -\frac{1-tq}{t^{\frac{1}{2}}(1-q)}\Bigg]^{n-1}(t^{\frac{1}{2}})^{wr(\widehat{\alpha })}\tau (\pi_{n}(\alpha ))$.
\end{center}
Here $wr(\widehat{\alpha })$ stands for the writhe number of the mixed link $\widehat{\alpha }$.

\end{defi}

\begin{lem}
HOMFLYPT polynomials of type $G(e,1)$ satisfy the following Skein relations:
\begin{center}
$\displaystyle (qt)^{-\frac{1}{2}}X_{L_{+}}-(qt)^{\frac{1}{2}}X_{L_{-}}=(q^{\frac{1}{2}}-q^{-\frac{1}{2}})X_{L_{0}}$,
\end{center}
\begin{center}
$X_{M_{e}}=a_{e-1}X_{M_{e-1}}+\cdots+a_{1}X_{M_{1}}+a_{0}X_{M_{0}}$.
\end{center}
Here, $a_{i}$ is defined by 
\begin{center}
$(T_{0}-u_{1})\cdots (T_{0}-u_{e})=T_{0}^{e}-a_{e-1}T_{0}^{e-1}-\cdots-a_{1}T_{0}-a_{0}$
\end{center}
and $L_{+},L_{-},L_{0},M_{e},\ldots,M_{1},M_{0}$ are given by the following local mixed link diagrams:
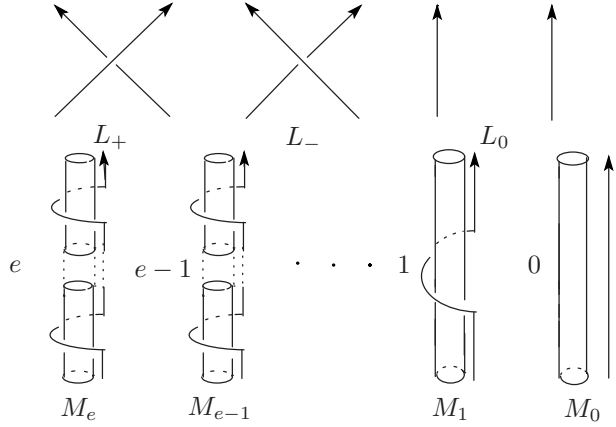
\begin{center}
\unitlength 0.1in
\begin{picture}( 26.2700,  6.4500)(  4.2000,-10.4500)
%
{\color[named]{Black}{%
\special{pn 8}%
\special{pa 420 990}%
\special{pa 1030 400}%
\special{fp}%
\special{sh 1}%
\special{pa 1030 400}%
\special{pa 968 432}%
\special{pa 992 438}%
\special{pa 996 462}%
\special{pa 1030 400}%
\special{fp}%
}}%
%
{\color[named]{Black}{%
\special{pn 8}%
\special{pa 1020 990}%
\special{pa 740 720}%
\special{fp}%
}}%
%
{\color[named]{Black}{%
\special{pn 8}%
\special{pa 690 680}%
\special{pa 420 410}%
\special{fp}%
\special{sh 1}%
\special{pa 420 410}%
\special{pa 454 472}%
\special{pa 458 448}%
\special{pa 482 444}%
\special{pa 420 410}%
\special{fp}%
}}%
%
{\color[named]{Black}{%
\special{pn 8}%
\special{pa 2010 1000}%
\special{pa 1420 400}%
\special{fp}%
\special{sh 1}%
\special{pa 1420 400}%
\special{pa 1452 462}%
\special{pa 1458 438}%
\special{pa 1482 434}%
\special{pa 1420 400}%
\special{fp}%
}}%
%
{\color[named]{Black}{%
\special{pn 8}%
\special{pa 1420 1000}%
\special{pa 1690 720}%
\special{fp}%
}}%
%
{\color[named]{Black}{%
\special{pn 8}%
\special{pa 1740 680}%
\special{pa 2020 410}%
\special{fp}%
\special{sh 1}%
\special{pa 2020 410}%
\special{pa 1958 442}%
\special{pa 1982 448}%
\special{pa 1986 472}%
\special{pa 2020 410}%
\special{fp}%
}}%
%
{\color[named]{Black}{%
\special{pn 8}%
\special{pa 3020 1010}%
\special{pa 3020 410}%
\special{fp}%
\special{sh 1}%
\special{pa 3020 410}%
\special{pa 3000 478}%
\special{pa 3020 464}%
\special{pa 3040 478}%
\special{pa 3020 410}%
\special{fp}%
}}%
%
{\color[named]{Black}{%
\special{pn 8}%
\special{pa 2420 990}%
\special{pa 2420 410}%
\special{fp}%
\special{sh 1}%
\special{pa 2420 410}%
\special{pa 2400 478}%
\special{pa 2420 464}%
\special{pa 2440 478}%
\special{pa 2420 410}%
\special{fp}%
\special{pa 2420 410}%
\special{pa 2420 410}%
\special{fp}%
\special{pa 2420 410}%
\special{pa 2420 410}%
\special{fp}%
}}%
\put(7.1000,-11.0000){\makebox(0,0){$L_{+}$}}%
\put(17.2000,-11.1000){\makebox(0,0){$L_{-}$}}%
\put(27.2000,-11.0000){\makebox(0,0){$L_{0}$}}%
\end{picture}%
\end{center}
\begin{center}
\unitlength 0.1in
\begin{picture}( 32.2200, 13.0000)( 13.9100,-17.6900)
%
{\color[named]{Black}{%
\special{pn 8}%
\special{pa 1888 1180}%
\special{pa 1888 1500}%
\special{fp}%
}}%
%
{\color[named]{Black}{%
\special{pn 8}%
\special{pa 1888 1564}%
\special{pa 1888 1650}%
\special{fp}%
\special{pa 1888 1650}%
\special{pa 1888 1650}%
\special{fp}%
}}%
%
{\color[named]{Black}{%
\special{pn 8}%
\special{pa 1736 1532}%
\special{pa 1736 1658}%
\special{fp}%
}}%
%
{\color[named]{Black}{%
\special{pn 8}%
\special{ar 1808 1658 80 32  6.0264770 6.2831853}%
\special{ar 1808 1658 80 32  0.0000000 3.1415927}%
}}%
%
{\color[named]{Black}{%
\special{pn 8}%
\special{ar 1808 1658 80 32  3.4217006 3.6379168}%
\special{ar 1808 1658 80 32  4.2865655 4.5027817}%
\special{ar 1808 1658 80 32  5.1514304 5.3676466}%
\special{ar 1808 1658 80 32  6.0162952 6.2325114}%
}}%
%
{\color[named]{Black}{%
\special{pn 8}%
\special{ar 1962 1438 300 78  1.6142472 3.6891549}%
}}%
%
{\color[named]{Black}{%
\special{pn 8}%
\special{ar 1950 1408 302 78  4.5510139 4.6737178}%
}}%
%
{\color[named]{Black}{%
\special{pn 8}%
\special{pa 1944 1338}%
\special{pa 1944 1188}%
\special{fp}%
}}%
%
{\color[named]{Black}{%
\special{pn 8}%
\special{pa 1936 1516}%
\special{pa 1936 1666}%
\special{fp}%
}}%
%
{\color[named]{Black}{%
\special{pn 8}%
\special{ar 1944 1416 264 86  3.8714923 3.9404578}%
\special{ar 1944 1416 264 86  4.1473544 4.2163199}%
\special{ar 1944 1416 264 86  4.4232164 4.4869634}%
}}%
%
{\color[named]{Black}{%
\special{pn 8}%
\special{pa 1736 1462}%
\special{pa 1736 1196}%
\special{fp}%
}}%
%
{\color[named]{Black}{%
\special{pn 8}%
\special{ar 1808 1180 74 24  1.0427219 6.2831853}%
\special{ar 1808 1180 74 24  0.0000000 0.7853982}%
}}%
%
{\color[named]{Black}{%
\special{pn 8}%
\special{pa 1896 512}%
\special{pa 1896 832}%
\special{fp}%
}}%
%
{\color[named]{Black}{%
\special{pn 8}%
\special{pa 1896 896}%
\special{pa 1896 982}%
\special{fp}%
\special{pa 1896 982}%
\special{pa 1896 982}%
\special{fp}%
}}%
%
{\color[named]{Black}{%
\special{pn 8}%
\special{pa 1744 864}%
\special{pa 1744 990}%
\special{fp}%
}}%
%
{\color[named]{Black}{%
\special{pn 8}%
\special{ar 1816 990 80 32  6.0264770 6.2831853}%
\special{ar 1816 990 80 32  0.0000000 3.1415927}%
}}%
%
{\color[named]{Black}{%
\special{pn 8}%
\special{ar 1816 990 80 32  3.4217006 3.6379168}%
\special{ar 1816 990 80 32  4.2865655 4.5027817}%
\special{ar 1816 990 80 32  5.1514304 5.3676466}%
\special{ar 1816 990 80 32  6.0162952 6.2325114}%
}}%
%
{\color[named]{Black}{%
\special{pn 8}%
\special{ar 1970 770 300 80  1.6142472 3.6711823}%
}}%
%
{\color[named]{Black}{%
\special{pn 8}%
\special{ar 1958 740 302 78  4.5562922 4.6766899}%
}}%
%
{\color[named]{Black}{%
\special{pn 8}%
\special{pa 1952 668}%
\special{pa 1952 520}%
\special{fp}%
}}%
%
{\color[named]{Black}{%
\special{pn 8}%
\special{pa 1946 848}%
\special{pa 1946 998}%
\special{fp}%
}}%
%
{\color[named]{Black}{%
\special{pn 8}%
\special{ar 1952 746 264 86  3.8685616 3.9371330}%
\special{ar 1952 746 264 86  4.1428473 4.2114187}%
\special{ar 1952 746 264 86  4.4171330 4.4853135}%
}}%
%
{\color[named]{Black}{%
\special{pn 8}%
\special{pa 1744 794}%
\special{pa 1744 528}%
\special{fp}%
}}%
%
{\color[named]{Black}{%
\special{pn 8}%
\special{ar 1816 512 74 24  1.0560122 6.2831853}%
\special{ar 1816 512 74 24  0.0000000 0.7853982}%
}}%
%
{\color[named]{Black}{%
\special{pn 8}%
\special{pa 1734 1002}%
\special{pa 1734 1234}%
\special{dt 0.045}%
\special{pa 1734 1234}%
\special{pa 1734 1234}%
\special{dt 0.045}%
\special{pa 1734 1234}%
\special{pa 1734 1234}%
\special{dt 0.045}%
}}%
%
{\color[named]{Black}{%
\special{pn 8}%
\special{pa 1896 986}%
\special{pa 1896 1182}%
\special{dt 0.045}%
\special{pa 1896 1182}%
\special{pa 1896 1182}%
\special{dt 0.045}%
}}%
%
{\color[named]{Black}{%
\special{pn 8}%
\special{pa 2616 1180}%
\special{pa 2616 1500}%
\special{fp}%
}}%
%
{\color[named]{Black}{%
\special{pn 8}%
\special{pa 2616 1564}%
\special{pa 2616 1650}%
\special{fp}%
\special{pa 2616 1650}%
\special{pa 2616 1650}%
\special{fp}%
}}%
%
{\color[named]{Black}{%
\special{pn 8}%
\special{pa 2464 1532}%
\special{pa 2464 1658}%
\special{fp}%
}}%
%
{\color[named]{Black}{%
\special{pn 8}%
\special{ar 2536 1658 80 32  6.0264770 6.2831853}%
\special{ar 2536 1658 80 32  0.0000000 3.1415927}%
}}%
%
{\color[named]{Black}{%
\special{pn 8}%
\special{ar 2536 1658 80 32  3.4217006 3.6379168}%
\special{ar 2536 1658 80 32  4.2865655 4.5027817}%
\special{ar 2536 1658 80 32  5.1514304 5.3676466}%
\special{ar 2536 1658 80 32  6.0162952 6.2325114}%
}}%
%
{\color[named]{Black}{%
\special{pn 8}%
\special{ar 2690 1438 300 78  1.6142472 3.6891549}%
}}%
%
{\color[named]{Black}{%
\special{pn 8}%
\special{ar 2678 1408 302 78  4.5542521 4.6766899}%
}}%
%
{\color[named]{Black}{%
\special{pn 8}%
\special{pa 2672 1338}%
\special{pa 2672 1188}%
\special{fp}%
}}%
%
{\color[named]{Black}{%
\special{pn 8}%
\special{pa 2666 1516}%
\special{pa 2666 1666}%
\special{fp}%
}}%
%
{\color[named]{Black}{%
\special{pn 8}%
\special{ar 2672 1416 264 86  3.8744078 3.9431757}%
\special{ar 2672 1416 264 86  4.1494794 4.2182473}%
\special{ar 2672 1416 264 86  4.4245510 4.4877795}%
}}%
%
{\color[named]{Black}{%
\special{pn 8}%
\special{pa 2464 1462}%
\special{pa 2464 1196}%
\special{fp}%
}}%
%
{\color[named]{Black}{%
\special{pn 8}%
\special{ar 2536 1180 74 24  1.0427219 6.2831853}%
\special{ar 2536 1180 74 24  0.0000000 0.7853982}%
}}%
%
{\color[named]{Black}{%
\special{pn 8}%
\special{pa 2624 512}%
\special{pa 2624 832}%
\special{fp}%
}}%
%
{\color[named]{Black}{%
\special{pn 8}%
\special{pa 2624 896}%
\special{pa 2624 982}%
\special{fp}%
\special{pa 2624 982}%
\special{pa 2624 982}%
\special{fp}%
}}%
%
{\color[named]{Black}{%
\special{pn 8}%
\special{pa 2472 864}%
\special{pa 2472 990}%
\special{fp}%
}}%
%
{\color[named]{Black}{%
\special{pn 8}%
\special{ar 2544 990 80 32  6.0264770 6.2831853}%
\special{ar 2544 990 80 32  0.0000000 3.1415927}%
}}%
%
{\color[named]{Black}{%
\special{pn 8}%
\special{ar 2544 990 80 32  3.4217006 3.6379168}%
\special{ar 2544 990 80 32  4.2865655 4.5027817}%
\special{ar 2544 990 80 32  5.1514304 5.3676466}%
\special{ar 2544 990 80 32  6.0162952 6.2325114}%
}}%
%
{\color[named]{Black}{%
\special{pn 8}%
\special{ar 2698 770 300 80  1.6141026 3.6711823}%
}}%
%
{\color[named]{Black}{%
\special{pn 8}%
\special{ar 2686 740 302 78  4.5562922 4.6766899}%
}}%
%
{\color[named]{Black}{%
\special{pn 8}%
\special{pa 2680 668}%
\special{pa 2680 520}%
\special{fp}%
}}%
%
{\color[named]{Black}{%
\special{pn 8}%
\special{pa 2674 848}%
\special{pa 2674 998}%
\special{fp}%
}}%
%
{\color[named]{Black}{%
\special{pn 8}%
\special{ar 2680 746 264 86  3.8685616 3.9371330}%
\special{ar 2680 746 264 86  4.1428473 4.2114187}%
\special{ar 2680 746 264 86  4.4171330 4.4853135}%
}}%
%
{\color[named]{Black}{%
\special{pn 8}%
\special{pa 2472 794}%
\special{pa 2472 528}%
\special{fp}%
}}%
%
{\color[named]{Black}{%
\special{pn 8}%
\special{ar 2544 512 74 24  1.0560122 6.2831853}%
\special{ar 2544 512 74 24  0.0000000 0.7853982}%
}}%
%
{\color[named]{Black}{%
\special{pn 8}%
\special{pa 2462 1002}%
\special{pa 2462 1234}%
\special{dt 0.045}%
\special{pa 2462 1234}%
\special{pa 2462 1234}%
\special{dt 0.045}%
\special{pa 2462 1234}%
\special{pa 2462 1234}%
\special{dt 0.045}%
}}%
%
{\color[named]{Black}{%
\special{pn 8}%
\special{pa 2624 986}%
\special{pa 2624 1182}%
\special{dt 0.045}%
\special{pa 2624 1182}%
\special{pa 2624 1182}%
\special{dt 0.045}%
}}%
%
{\color[named]{Black}{%
\special{pn 8}%
\special{pa 3832 510}%
\special{pa 3832 1254}%
\special{fp}%
}}%
%
{\color[named]{Black}{%
\special{pn 8}%
\special{pa 3832 1430}%
\special{pa 3832 1636}%
\special{fp}%
\special{pa 3832 1636}%
\special{pa 3832 1636}%
\special{fp}%
}}%
%
{\color[named]{Black}{%
\special{pn 8}%
\special{pa 3678 1356}%
\special{pa 3678 1656}%
\special{fp}%
}}%
%
{\color[named]{Black}{%
\special{pn 8}%
\special{ar 3906 1130 300 188  1.6141026 3.6903058}%
}}%
%
{\color[named]{Black}{%
\special{pn 8}%
\special{ar 3894 1080 302 188  4.5504853 4.6733698}%
}}%
%
{\color[named]{Black}{%
\special{pn 8}%
\special{pa 3880 1318}%
\special{pa 3880 1676}%
\special{fp}%
}}%
%
{\color[named]{Black}{%
\special{pn 8}%
\special{ar 3884 1096 266 206  3.8688846 3.9198400}%
\special{ar 3884 1096 266 206  4.0727062 4.1236616}%
\special{ar 3884 1096 266 206  4.2765279 4.3274833}%
\special{ar 3884 1096 266 206  4.4803495 4.4861414}%
}}%
%
{\color[named]{Black}{%
\special{pn 8}%
\special{pa 3678 1186}%
\special{pa 3678 550}%
\special{fp}%
}}%
%
{\color[named]{Black}{%
\special{pn 8}%
\special{ar 3758 1646 72 44  5.6135464 6.2831853}%
\special{ar 3758 1646 72 44  0.0000000 3.1415927}%
}}%
%
{\color[named]{Black}{%
\special{pn 8}%
\special{ar 3758 1646 72 44  3.3007294 3.5094250}%
\special{ar 3758 1646 72 44  4.1355120 4.3442076}%
\special{ar 3758 1646 72 44  4.9702946 5.1789902}%
\special{ar 3758 1646 72 44  5.8050772 6.0137729}%
}}%
%
{\color[named]{Black}{%
\special{pn 8}%
\special{pa 3678 1492}%
\special{pa 3678 1268}%
\special{fp}%
}}%
%
{\color[named]{Black}{%
\special{pn 8}%
\special{pa 3830 1544}%
\special{pa 3830 1338}%
\special{fp}%
}}%
%
{\color[named]{Black}{%
\special{pn 8}%
\special{ar 3750 504 82 36  0.0000000 6.2831853}%
}}%
%
{\color[named]{Black}{%
\special{pn 8}%
\special{pa 3678 1148}%
\special{pa 3678 1216}%
\special{fp}%
}}%
%
{\color[named]{Black}{%
\special{pn 8}%
\special{pa 3678 608}%
\special{pa 3678 522}%
\special{fp}%
}}%
%
{\color[named]{Black}{%
\special{pn 8}%
\special{pa 3884 900}%
\special{pa 3884 488}%
\special{fp}%
\special{sh 1}%
\special{pa 3884 488}%
\special{pa 3864 554}%
\special{pa 3884 540}%
\special{pa 3904 554}%
\special{pa 3884 488}%
\special{fp}%
}}%
%
{\color[named]{Black}{%
\special{pn 8}%
\special{pa 2678 556}%
\special{pa 2678 480}%
\special{fp}%
\special{sh 1}%
\special{pa 2678 480}%
\special{pa 2658 546}%
\special{pa 2678 532}%
\special{pa 2698 546}%
\special{pa 2678 480}%
\special{fp}%
}}%
%
{\color[named]{Black}{%
\special{pn 8}%
\special{pa 1950 624}%
\special{pa 1940 480}%
\special{fp}%
\special{sh 1}%
\special{pa 1940 480}%
\special{pa 1924 548}%
\special{pa 1944 532}%
\special{pa 1964 544}%
\special{pa 1940 480}%
\special{fp}%
}}%
%
{\color[named]{Black}{%
\special{pn 8}%
\special{pa 2670 980}%
\special{pa 2670 1186}%
\special{dt 0.045}%
}}%
%
{\color[named]{Black}{%
\special{pn 8}%
\special{pa 1940 960}%
\special{pa 1940 1208}%
\special{dt 0.045}%
}}%
%
{\color[named]{Black}{%
\special{pn 8}%
\special{pa 4326 1366}%
\special{pa 4326 1664}%
\special{fp}%
}}%
%
{\color[named]{Black}{%
\special{pn 8}%
\special{pa 4326 1196}%
\special{pa 4326 558}%
\special{fp}%
}}%
%
{\color[named]{Black}{%
\special{pn 8}%
\special{ar 4406 1654 72 44  5.6135464 6.2831853}%
\special{ar 4406 1654 72 44  0.0000000 3.1415927}%
}}%
%
{\color[named]{Black}{%
\special{pn 8}%
\special{ar 4406 1654 72 44  3.3247035 3.5333991}%
\special{ar 4406 1654 72 44  4.1594861 4.3681817}%
\special{ar 4406 1654 72 44  4.9942687 5.2029643}%
\special{ar 4406 1654 72 44  5.8290513 6.0377469}%
}}%
%
{\color[named]{Black}{%
\special{pn 8}%
\special{pa 4326 1500}%
\special{pa 4326 1278}%
\special{fp}%
}}%
%
{\color[named]{Black}{%
\special{pn 8}%
\special{ar 4398 514 82 34  0.0000000 6.2831853}%
}}%
%
{\color[named]{Black}{%
\special{pn 8}%
\special{pa 4326 1156}%
\special{pa 4326 1226}%
\special{fp}%
}}%
%
{\color[named]{Black}{%
\special{pn 8}%
\special{pa 4326 616}%
\special{pa 4326 530}%
\special{fp}%
}}%
%
{\color[named]{Black}{%
\special{pn 8}%
\special{pa 4326 856}%
\special{pa 4326 1346}%
\special{fp}%
}}%
%
{\color[named]{Black}{%
\special{pn 8}%
\special{pa 4586 1664}%
\special{pa 4586 504}%
\special{fp}%
\special{sh 1}%
\special{pa 4586 504}%
\special{pa 4566 572}%
\special{pa 4586 558}%
\special{pa 4606 572}%
\special{pa 4586 504}%
\special{fp}%
}}%
%
{\color[named]{Black}{%
\special{pn 4}%
\special{sh 1}%
\special{ar 2958 1062 6 6 0  6.28318530717959E+0000}%
\special{sh 1}%
\special{ar 3146 1072 6 6 0  6.28318530717959E+0000}%
\special{sh 1}%
\special{ar 3336 1072 6 6 0  6.28318530717959E+0000}%
\special{sh 1}%
\special{ar 3336 1072 6 6 0  6.28318530717959E+0000}%
}}%
%
{\color[named]{Black}{%
\special{pn 4}%
\special{sh 1}%
\special{ar 3336 1072 6 6 0  6.28318530717959E+0000}%
\special{sh 1}%
\special{ar 3336 1072 6 6 0  6.28318530717959E+0000}%
}}%
\put(14.8100,-10.7900){\makebox(0,0){$e$}}%
\put(22.5500,-10.9600){\makebox(0,0){$e-1$}}%
\put(35.1500,-10.7100){\makebox(0,0){$1$}}%
\put(41.9900,-10.7100){\makebox(0,0){$0$}}%
\put(18.1400,-18.2600){\makebox(0,0){$M_{e}$}}%
\put(25.5200,-18.2600){\makebox(0,0){$M_{e-1}$}}%
\put(37.5800,-18.2600){\makebox(0,0){$M_{1}$}}%
\put(44.5100,-18.3400){\makebox(0,0){$M_{0}$}}%
%
{\color[named]{Black}{%
\special{pn 8}%
\special{pa 4470 1640}%
\special{pa 4470 510}%
\special{fp}%
}}%
\end{picture}%

\end{center}
\end{lem}

\section{Alexander polynomials for mixed links}
In this section we define the Alexander polynomials for mixed links.

If we want to define the Alexander polynomial for a mixed link, we must consider the specialization $t\rightarrow q^{-1}$.
However, a priori, the specialization does not make sense since the term $1-qt$ appears in $X(q,t)$.
To solve this problem we can use the Skein relations for $X(q,t)$. As explained in the previous section the link polynomials $X(q,t)$ satisfy the Skein relation. Conversely by giving some initial conditions, we can define $X(q,t)$ by using the Skein relations. In particular, we can
define the Alexander polynomial for a mixed link. However, this definition is indirect and is not explicit. 

In fact, thanks to the result of Geck-Iancu-Malle explained in Section 3, we can define the Alexander polynomial for a mixed link explicitly.

To see this let us focus on the term $R^{\bm{\lambda}}$. 

After the change of variables
\begin{center}
$\displaystyle z=-\frac{1-q}{1-tq}$,
\end{center}
$R^{\bm{\lambda}}$ can be written as
$\widehat{g}^{\bm{\lambda}}/(1-tq)^{n}$ where $\widehat{g}^{\bm{\lambda}}$ is given by
\begin{center}
$\displaystyle \prod_{k=1}^{e}\prod_{x\in \lambda _{k}}\Bigg( -(1-q)(1-q^{c(x)})\prod_{\stackrel{\scriptstyle 1\leq l\leq e}{l\not= k}}(q^{c(x)}u_{k}-u_{l})+(1-tq)(1-q)C(e,k,x)\Bigg)$.
\end{center} 
Since at least one component of $\bm{\lambda}$,  say $\lambda _{k}$, is not empty, when $x$ is the (1,1)-component of $\lambda _{k}$,
the corresponding term is given by $(1-tq)(1-q)C(e,k,x)$. This allows us to evaluate at $t=q^{-1}$.
\begin{defi}
Let $\widehat{\alpha} $ be a mixed link obtained as the closure of a mixed braid $\alpha $. We define \textbf{Alexander polynomial of type $\bm{G(e,1)}$} for $\widehat{\alpha }$ by
\begin{center}
 $\Delta^{G(e,1)}(\widehat{\alpha })=X_{\widehat{\alpha }}(q,q^{-1})$.
\end{center}
\end{defi}

The following is a direct consequence of Lemma 4.7.
\begin{cor}
Alexander polynomials of type $G(e,1)$ satisfy the following Skein relations:
\begin{center}
$\displaystyle \Delta^{G(e,1)}_{L_{+}}-\Delta^{G(e,1)}_{L_{-}}=(q^{\frac{1}{2}}-q^{-\frac{1}{2}})\Delta^{G(e,1)}_{L_{0}}$,
\end{center}
\begin{center}
$\Delta^{G(e,1)}_{M_{e}}=a_{e-1}\Delta^{G(e,1)}_{M_{e-1}}+\cdots+a_{1}\Delta^{G(e,1)}_{M_{1}}+a_{0}\Delta^{G(e,1)}_{M_{0}}$.
\end{center}
Here, the local mixed link diagrams $L_{+},L_{-},L_{0},M_{e},M_{e-1},\ldots,M_{1},M_{0}$ are as in Lemma 4.7.
\end{cor}

Finally we discuss simplification of the Alexander polynomials of type $G(e,1)$.

Assume that $\bm{\lambda}$ has at least two non-empty components. Then the same consideration shows that after the specialization $t=q^{-1}$,
the corresponding $R^{\bm{\lambda}}$ is zero. Similarly, if $\bm{\lambda}$ has a component which has at least two diagonal boxes, the corresponding
$R^{\bm{\lambda}}$ is also zero. This shows that when we consider the Alexander polynomial of type $G(e,1)$, we only have to consider the
$e$-Young diagrams which have the form $(\emptyset,\cdots,\emptyset,\lambda^{(a)},\emptyset,\cdots,\emptyset)$, where $\lambda ^{(a)}$ is
the Young diagram $(a+1,1,\cdots,1)$ of size $n$.

Summarizing the above argument we have the following lemma.
\begin{lem}
Let $\alpha $ be a mixed braid. Then we have \\
\linebreak
$\Delta^{G(e,1)}(\widehat{\alpha })$\\
$\displaystyle =(q-1)^{-(n-1)}q^{-\frac{n-wr(\widehat{\alpha})-1}{2}}\sum_{\stackrel{\scriptstyle 1\leq p\leq e}{0\leq a\leq n-1}}(-1)^{en}
\Big( \prod_{\stackrel{\scriptstyle 1\leq i\leq e}{i\not= p}}
u_{i}^{-n}\Big) D^{\bm{\lambda}_{p}^{(a)}}g^{ \bm{\lambda}_{p}^{(a)}}\chi _{\bm{\lambda}_{p}^{(a)}}$.\\
\linebreak
Here 
\[
\bm{\lambda}_{p}^{(a)}=(\emptyset,\cdots,\emptyset,\stackrel{p}{\breve{\lambda}^{(a)}},\emptyset,\cdots,\emptyset ),
\]
\[
g^{ \bm{\lambda}_{p}^{(a)}} =(-1)^{n-1}(1-q)C(e,p)\prod_{\stackrel{\scriptstyle x\in \lambda ^{(a)}}{x\not= (1,1)}}\Bigg( (1-q)(1-q^{c(x)})
\prod_{\stackrel{\scriptstyle 1\leq l\leq e}{l\not= p}}(q^{c(x)}u_{k}-u_{l})\Bigg),\\
\]
\[
C(e,p)=\prod_{i=1}^{e-1}(u_{k}^{e-i-1}\sum_{j=1}^{i}(-1)^{i-j}y_{j}\sigma _{i-j})+(-1)^{e-1}\prod_{\stackrel{\scriptstyle 1\leq l\leq e}{l\not= k}}u_{l}.
\]
\end{lem}

\section{Quantum calculus}
In this section we calculate $D^{\bm{\lambda}_{p}^{(a)}}g^{ \bm{\lambda}_{p}^{(a)}}$.  

We first consider the case of $p\in \{ 2,\cdots,e \}$. In this case $\alpha _{i,j}$ are given by
\begin{center}
$\alpha_{1,i}=n-i+1\ \ (1\leq i\leq n+1),\ \ \ \ \alpha _{l,i}=n-i\ \ (l\not=1,p,\ 1\leq i\leq n)$,
\begin{eqnarray}
\alpha_{p,i}= \left\{ \begin{array}{ll}
a+n & (i=1) \\
n-i+1 & (2 \leq j\leq b+1) \\
n-i  &  (b+2\leq n).\\
\end{array} \right. \nonumber
\end{eqnarray} 
\end{center}
Here $b=n-a-1$.

\begin{lem}
For $l\not= 1,p $ and $i,j\geq 2$, we have the following.
\begin{itemize}
\item[(1)] $\displaystyle \prod_{j'=1}^{j-1}(q^{\alpha _{1,j'}}-q^{\alpha _{1,j}})=q^{(j-1)(n-j+1)}\prod_{h=1}^{j-1}(q^{h}-1).$
\item[(2)] $\displaystyle \prod_{j'=1}^{j-1}(q^{\alpha _{l,j'}}-q^{\alpha _{l,j}})=q^{(j-1)(n-j)}\prod_{h=1}^{j-1}(q^{h}-1).$
\item[(3)] $\displaystyle \prod_{i'=1}^{i-1}(q^{\alpha _{p,i'}}-q^{\alpha _{p,i}})$\\
$= \left\{ \begin{array}{ll}
\displaystyle q^{(i-1)(n-i)}\frac{q^{a+i}-1}{q^{i-b-1}-1} \prod_{h=1}^{i-1}(q^{h}-1) & (b+3\leq i\leq n), \\
\displaystyle q^{(i-1)(n-i)}\frac{q^{a+i}-1}{q-1} \prod_{h=1}^{i-1}(q^{h}-1) & (i=b+2), \\
\displaystyle q^{(i-1)(n-i)+i-1}\frac{q^{a+i-1}-1}{q^{i-1}-1} \prod_{h=1}^{i-1}(q^{h}-1)  &  (2\leq i\leq b+1).\\
\end{array} \right. \nonumber $
\end{itemize}
\begin{proof}
We only prove the 1st identity in (3).  
\begin{align*}
(\text{LHS}) & =\prod_{i'=1}^{i-1}(q^{\alpha _{p,i'}}-q^{n-i}) \\  
             & =q^{(i-1)(n-i)}\prod_{i'=1}^{i-1}(q^{\alpha_{p,i'}-(n-i)}-1)  \\ 
             & =q^{(i-1)(n-i)}(q^{a+i}-1)\prod_{i'=2}^{b+1}(q^{(n-i'+1)-(n-i)}-1)\prod_{i'=b+2}^{i-1}(q^{(n-i')-(n-i)}-1) \\ 
             & =q^{(i-1)(n-i)}(q^{a+i}-1)\prod_{i'=2}^{b+1}(q^{i-i'+1}-1)\prod_{i'=b+2}^{i-1}(q^{i-i'}-1) \\ 
             & =q^{(i-1)(n-i)}\frac{q^{a+i}-1}{q^{i-b-1}-1} \prod_{h=1}^{i-1}(q^{h}-1). 
\end{align*}
\end{proof}
\end{lem}
Now we can decompose  $D^{\bm{\lambda}_{p}^{(a)}}$ into the following three factors:\\
\linebreak
(1)\ \ \ $\displaystyle \frac{\displaystyle (-1)^{\binom{e}{2}\binom{n}{2}+n(e-1)} q^{-n}\prod_{1\leq l\leq e}u_{l}^{n}}
{\displaystyle q^{f(n,e)}\prod_{2\leq l\leq e}(q^{-1}u_{1}-u_{l})^{n}\prod_{2\leq k<l\leq e}(u_{k}-u_{l})^{n}}.$\\
\linebreak
\linebreak
(2) $\displaystyle  \prod_{\stackrel{\scriptstyle \alpha ,\alpha '\in A_{1}}{\alpha >\alpha'}}(q^{\alpha }-q^{\alpha'})\times q^{-\binom{n+1}{2}}u_{1}^{\binom{n+1}{2}}$
$\hfill(k=l=1),$\\
$\displaystyle \prod_{\stackrel{\scriptstyle 2\leq l\leq e}{l\not =p}}\prod_{(\alpha ,\alpha ')\in \bold{A}_{1}\times \bold{A}_{l}}(q^{\alpha -1}u_{1}-q^{\alpha '}u_{l})
$ $\hfill (k=1,l\geq 2,l\not=p),$\\
$\displaystyle \prod_{\stackrel{\scriptstyle 2\leq l\leq e}{l\not =p}}\prod_{\alpha ,\alpha '\in \bold{A}_{l},\alpha >\alpha '}(q^{\alpha}-q^{\alpha '})
\times \prod_{\stackrel{\scriptstyle 2\leq l\leq e}{l\not= p}}u_{l}^{n}$ 
$\hfill (k=l\geq 2,l\not= p),$\\ 
$\displaystyle \prod_{\stackrel{\scriptstyle 2\leq k<l\leq e}{k,l\not= p}}\prod_{(\alpha ,\alpha ')\in \bold{A}_{k}\times \bold{A}_{l}}(q^{\alpha} u_{k}-q^{\alpha'}u_{l})$
$\hfill (2\leq k<l\leq e,k,l\not=p),$\\
$\displaystyle \prod_{(\alpha ,\alpha ')\in \bold{A}_{1}\times \bold{A}_{p}}(q^{\alpha -1}u_{1}-q^{\alpha '}u_{p})$
$\hfill (k=1,l=p),$\\
$\displaystyle \prod_{\stackrel{\scriptstyle \alpha ,\alpha '\in \bold{A}_{p}}{\alpha >\alpha '}}(q^{\alpha }-q^{\alpha '})\times u_{p}^{\binom{n}{2}}$
$\hfill (k=l=p),$\\
$\displaystyle \prod_{\stackrel{\scriptstyle 2\leq k<l\leq e}{k=p\ \text{or}\ l=p}}\prod_{(\alpha ,\alpha ')\in \bold{A}_{k}\times \bold{A}_{l}}(q^{\alpha }u_{k}-q^{\alpha '}
u_{l})$ 
$\hfill (2\leq k<l\leq e;k=p\ \text{or}\ l=p).$\\
\linebreak
\linebreak
(3) $\displaystyle \prod_{\alpha \in \bold{A}_{1}}\prod_{h=1}^{\alpha }(q^{h}-1)\times q^{-\sum_{\alpha \in \bold{A}_{1}}\alpha }u_{1}^{\sum_{\alpha \in \bold{A}_{1}}\alpha}$
$\hfill (k=l=p),$\\
$\displaystyle \prod_{2\leq l\leq e}\prod_{\alpha \in \bold{A}_{1}}\prod_{h=1}^{\alpha }(q^{h-1}u_{1}-u_{l})$
$\hfill (k=1,l\geq 2),$\\
$\displaystyle \prod_{\stackrel{\scriptstyle 2\leq l\leq e}{l\not= p}}\prod_{\alpha \in \bold{A}_{l}}\prod_{h=1}^{\alpha }(q^{h}-1)\times \prod_{2\leq l\leq e,l\not= p}u_{l}^{\sum_{\alpha \in \bold{A}_{l}}\alpha }$
$\hfill (k=l\geq 2,l\not= p),$\\
$\displaystyle \prod_{\stackrel{\scriptstyle 2\leq k\leq e}{k\not= p}}\prod_{\alpha \in \bold{A}_{k}}\prod_{h=1}^{\alpha }(q^{h}u_{k}-u_{1})$  
$\hfill (k\geq 2,k\not= p,l=1),$\\
$\displaystyle \prod_{\stackrel{\scriptstyle k,l\geq 2}{k\not =p,l}}\prod_{\alpha \in \bold{A}_{k}}\prod_{h=1}^{\alpha }(q^{h}u_{k}-u_{l})$
$\hfill (k\not= l,p;k,l\geq 2),$\\ 
$\displaystyle \prod_{\alpha \in \bold{A}_{p}}\prod_{h=1}^{\alpha }(q^{h}u_{p}-u_{1})$
$\hfill (k=p,l=1),$\\
$\displaystyle \prod_{\alpha \in \bold{A}_{p}}\prod_{h=1}^{\alpha }(q^{h}-1)\times u_{p}^{ \sum_{\alpha\in \bold{A}_{p}}\alpha }$
$\hfill (k=l=p),$\\
$\displaystyle \prod_{\stackrel{\scriptstyle 2\leq l\leq e}{l\not= p}}\prod_{\alpha \in \bold{A}_{p}}\prod_{h=1}^{\alpha }(q^{h}u_{p}-u_{l})$
$\hfill (k=p,l\not= p,l\geq 2).$\\
\linebreak 
We combine the terms which relate $\bold{A}_{p}$, that is, we put 
\begin{center}
$\displaystyle D_{0}^{\bm{\lambda}_{p}^{(a)}}=\frac{\displaystyle \prod_{\alpha ,\alpha '\in \bold{A}_{p}}(q^{\alpha }-q^{\alpha '})}{\displaystyle
\prod_{\alpha \in \bold{A}_{p}}\prod_{h=1}^{\alpha }(q^{h}-1)},$
\end{center}
$\displaystyle D_{1}^{\bm{\lambda}_{p}^{(a)}}= \frac{\displaystyle \prod_{(\alpha ,\alpha ')\in \bold{A}_{1}\times \bold{A}_{p}}(q^{\alpha -1}u_{1}-q^{\alpha '}u_{p})
\prod_{\stackrel{\scriptstyle 2\leq k<l\leq e}{k=p\ \text{or}\ l=p }} \prod_{(\alpha ,\alpha ')\in \bold{A}_{k}\times \bold{A}_{l}}(q^{\alpha }u_{k}-q^{\alpha '}u_{l})}
{\displaystyle \prod_{\alpha \in \bold{A}_{p}}\prod_{h=1}^{\alpha }(q^{h}u_{p}-q^{-1}u_{1})
\prod_{\stackrel{\scriptstyle 2\leq l\leq e}{l\not=p}}\prod_{\alpha \in \bold{A}_{p}}\prod_{h=1}^{\alpha }(q^{h}u_{p}-u_{l})},$
\begin{center}
$\displaystyle g_{0}^{\bm{\lambda}_{p}^{(a)}}=(-1)^{n-1}(1-q)^{n}C(e,p)\prod_{i=1}^{a}(1-q^{i})\prod_{j=1}^{b}(1-q^{-j}),$\\
$\displaystyle g_{1}^{\bm{\lambda}_{p}^{(a)}}=\prod_{i=1}^{a}\prod_{\stackrel{\scriptstyle 1\leq l\leq e}{l\not= p}}(q^{i}u_{p}-u_{l})
\prod_{j=1}^{b}\prod_{\stackrel{\scriptstyle 1\leq l\leq e}{l\not= p}}(q^{-j}u_{p}-u_{l}).$
\end{center}
We note that $g^{\bm{\lambda}_{p}^{(a)}}=g_{0}^{\bm{\lambda}_{p}^{(a)}}g_{1}^{\bm{\lambda}_{p}^{(a)}}$.
\begin{lem}We have
\begin{center}
$\displaystyle D_{0}^{\bm{\lambda}_{p}^{(a)}}g_{0}^{\bm{\lambda}_{p}^{(a)}}=(-1)^{a+n}C(e,p)q^{\frac{n(n-1)(n-2)}{6}}\frac{(1-q)^{n}}{1-q^{n}}$.
\end{center}
\begin{proof}
By Lemma 6.1 (3) we have\\
\linebreak
$\displaystyle \prod_{\alpha ,\alpha '\in \bold{A}_{p}}(q^{\alpha }-q^{\alpha '}) $\\
$\displaystyle =\prod_{i=2}^{n}\prod_{i'=1}^{i}(q^{\alpha _{p,i'}}-q^{\alpha _{p,i}})$ \\
$\displaystyle =q^{\sum_{i=2}^{n}(n-i)(i-1)+\sum_{i=2}^{b+1}(i-1)}\prod_{i=2}^{n}\prod_{h=1}^{i-1}(q^{h}-1)$\\
$\hspace{20mm}$ $\displaystyle \times \prod_{i=2}^{b+1}\frac{q^{a+i-1}-1}{q^{i-1}-1}\times \frac{q^{n+1}-1}{q-1}\times \prod_{i=b+3}^{n}\frac{q^{a+i}-1}{q^{i-b-1}-1}.$\\
\linebreak                                                                     
On the other hand, since
\begin{center}
$\displaystyle \prod_{j=1}^{b}(1-q^{-j})=q^{-\sum_{j=1}^{b}j}\prod_{j=1}^{b}(q^{j}-1)$,\ \ \ \ $\displaystyle \prod_{i=1}^{a}(1-q^{i})=(-1)^{a}\prod_{i=1}^{a}(q^{i}-1)$,
\end{center} 
we have 
\begin{center}
$\displaystyle \prod_{i=2}^{b+1}\frac{q^{a+i-1}-1}{q^{i-1}-1}\prod_{j=1}^{b}(1-q^{-j})=q^{-\sum_{j=1}^{b}j}\prod_{i=a+1}^{n-1}(q^{i}-1)$,\\
$\displaystyle \prod_{i=b+3}^{n}\frac{q^{a+i}-1}{q^{i-b-1}-1}\prod_{i=1}^{a}(1-q^{i})=(-1)^{a}(q-1)\prod_{i=n+2}^{a+n}(q^{i}-1)$.
\end{center}
By combining these calculations, we also have\\
\linebreak
$D_{0}^{\bm{\lambda}_{p}^{(a)}}g_{0}^{\bm{\lambda}_{p}^{(a)}}$\\
$\displaystyle =(-1)^{a+n-1}C(e,p)q^{\sum_{i=2}^{n}(n-i)(i-1)}\frac{(1-q)^{n}}{q^{n}-1}
\frac{\displaystyle \prod_{i=2}^{n}\prod_{h=1}^{i-1}(q^{h}-1)\prod_{i=a+1}^{a+n}(q^{i}-1)}{ \displaystyle \prod_{\alpha \in \bold{A}_{p}}\prod_{h=1}^{\alpha }(q^{h}-1)}$.\\
\linebreak
Since
\begin{align*}
\prod_{\alpha \in \bold{A}_{p}}\prod_{h=1}^{\alpha }(q^{h}-1) & = \frac{\displaystyle\prod_{i=1}^{n}\prod_{h=1}^{i-1}(q^{h}-1)\prod_{i=a+1}^{a+n}(q^{i}-1)}
{\displaystyle\prod_{h=1}^{n-b-1}(q^{h}-1)} \\
& = \prod_{i=2}^{n}\prod_{h=1}^{i-1}(q^{h}-1)\prod_{i=a+1}^{a+n}(q^{i}-1),
\end{align*}
we obtain the desired formula.
\end{proof}
\end{lem}
\begin{lem}We have\\
\linebreak
$D_{1}^{\bm{\lambda}_{p}^{(a)}}g_{1}^{\bm{\lambda}_{p}^{(a)}}$\\
$\displaystyle =(-1)^{np}\frac{1}{u_{p}-u_{1}}\prod_{j=1}^{n}\prod_{h=j-n}^{1}(q^{h}u_{p}-u_{1})
\prod_{\stackrel{\scriptstyle 2\leq l\leq e}{l\not=p}}\frac{1}{u_{p}-u_{l}}\prod_{i=1}^{n}\prod_{h=i-n}^{0}(q^{h}u_{p}-u_{l}).$
\begin{proof}
We first focus on the factors of $D_{1}^{\bm{\lambda}_{p}^{(a)}}$ which have the form $q^{*}u_{1}-q^{*}u_{p}$. 
It is given by
\begin{center}
$\displaystyle \frac{\displaystyle \prod_{(\alpha ,\alpha ')\in \bold{A}_{1}\times \bold{A}_{p}}(q^{\alpha -1}u_{1}-q^{\alpha '}u_{p})}
{\displaystyle \prod_{\alpha \in \bold{A}_{p}}\prod_{h=1}^{\alpha }(q^{h}u_{p}-q^{-1}u_{1})}$.
\end{center}
Now the numerator is equal to\\
\linebreak
$\displaystyle \prod_{i=1}^{n+1}\prod_{j=1}^{n}(q^{-1+(i-1)}u_{1}-q^{\alpha _{p,j}}u_{p})$\\
$\displaystyle =q^{n\sum_{i=1}^{n+1}(i-2)}\prod_{i=1}^{n+1}\prod_{j=1}^{n}(q^{\alpha _{p,j}}u_{p}-u_{1})$\\
$\displaystyle =q^{\frac{n(n+1)(n-2)}{2}}\frac{\displaystyle \prod_{i=1}^{n+1}\prod_{j=1}^{n}(q^{(j-1)-(i-2)}u_{p}-u_{1})\prod_{i=1}^{n+1}(q^{(a+n)-(i-2)}u_{p}-u_{1})}
{\displaystyle \prod_{i=1}^{n+1}(q^{(n-b-1)-(i-2)}u_{p}-u_{1})}$\\
$\displaystyle =q^{\frac{n(n+1)(n-2)}{2}}\frac{\displaystyle \prod_{i=1}^{n+1}\prod_{j=1}^{n}(q^{j-i+1}u_{p}-u_{1})\prod_{i=1}^{n+1}(q^{a+n-i+2}u_{p}-u_{1})}
{\displaystyle \prod_{i=1}^{n+1}(q^{n-i-b+1}u_{p}-u_{1})}$\\
\linebreak
and the denominator is equal to
\begin{align*}
\prod_{\alpha \in \bold{A}_{p}}\prod_{h=1}^{\alpha }(q^{h}u_{p}-q^{-1}u_{1}) & = \frac{\displaystyle \prod_{j=1}^{n}\prod_{h=1}^{j-1}(q^{h}u_{p}-q^{-1}u_{1})
\prod_{h=1}^{a+n}(q^{h}u_{p}-q^{-1}u_{1})}{\displaystyle \prod_{i=1}^{n-b-1}(q^{h}u_{p}-q^{-1}u_{1})}\\
& =q^{-\frac{n(n+1)}{2}}\frac{\displaystyle \prod_{j=1}^{n}\prod_{h=1}^{j-1}(q^{h+1}u_{p}-u_{1})
\prod_{h=1}^{a+n}(q^{h+1}u_{p}-u_{1})}{\displaystyle \prod_{i=1}^{n-b-1}(q^{h+1}u_{p}-u_{1})}.
\end{align*}
Thus the factor is given by
\begin{center}
$\displaystyle 
q^{\frac{n(n+1)(n-1)}{2}}\frac{\displaystyle \prod_{i=1}^{n+1}\prod_{j=1}^{n}(q^{j-i+1}u_{p}-u_{1})\prod_{i=1}^{n+1}(q^{a+n-i+2}u_{p}-u_{1})\prod_{i=1}^{n-b-1}(q^{h+1}u_{p}-u_{1})}
{\displaystyle \prod_{i=1}^{n+1}(q^{n-i-b+1}u_{p}-u_{1}) \prod_{j=1}^{n}\prod_{h=1}^{j-1}(q^{h+1}u_{p}-u_{1})
\prod_{h=1}^{a+n}(q^{h+1}u_{p}-u_{1})}$.
\end{center}
Now the following three identities hold:
\begin{center}
$\displaystyle \frac{\displaystyle \prod_{i=1}^{n+1}\prod_{j=1}^{n}(q^{j-i+1}u_{p}-u_{1})}{\displaystyle\prod_{h=1}^{j-1}(q^{h+1}u_{p}-u_{1})}=\prod_{j=1}^{n}\prod_{h=j-n}^{1}
(q^{h}u_{p}-u_{1}),$\\
$\displaystyle \frac{\displaystyle\prod_{i=1}^{n+1}(q^{a+n-i+2}u_{p}-u_{1})}{\displaystyle\prod_{h=1}^{a+n}(q^{h+1}u_{p}-u_{1})}=(qu_{p}-u_{1})\Bigg[ \prod_{i=1}^{a}(q^{i}u_{p}-u_{1})\Bigg]^{-1},$\\
$\displaystyle \frac{\displaystyle\prod_{i=1}^{n-b-1}(q^{h+1}u_{p}-u_{1})}{\displaystyle\prod_{i=1}^{n+1}(q^{n-i-b+1}u_{p}-u_{1})}=\frac{1}{(u_{p}-u_{1})(qu_{p}-u_{1})}
\Bigg[\prod_{j=1}^{b}(q^{-j}u_{p}-u_{1})\Bigg]^{-1}$.
\end{center}
Thus by combining with the corresponding factors in $g_{1}^{\bm{\lambda}_{p}^{(a)}}$, we have

\[
q^{n(n+1)(n-1)}\frac{1}{u_{p}-u_{1}}\prod_{j=1}^{n}\prod_{h=j-n}^{1}(q^{h}u_{p}-u_{1}).
\]

Next we focus on the factors in $D_{1}^{\bm{\lambda}_{p}^{(a)}}$ which have the form $q^{*}u_{p}-q^{*}u_{l}\ (l\not= 1,p)$.
It is given by
\[ 
\frac{\displaystyle\prod_{\stackrel{\scriptstyle 2\leq k<l\leq e}{k=p\ \text{or}\ l=p}}\prod_{(\alpha ,\alpha ')\in \bold{A}_{k}\times \bold{A}_{l}}(q^{\alpha }u_{k}-q^{\alpha '}u_{l})}
{\displaystyle\prod_{\stackrel{\scriptstyle 2\leq l\leq e}{l\not= p}}\prod_{\alpha \in \bold{A}_{p}}\prod_{h=1}^{\alpha }(q^{h}u_{p}-u_{l})}.
\]
We note that the numerator is equal to
\[
(-1)^{np}\prod_{\stackrel{\scriptstyle 2\leq l\leq e}{l\not= p}}\prod_{(\alpha ,\alpha ')\in \bold{A}_{p}\times \bold{A}_{l}}(q^{\alpha }u_{p}-q^{\alpha '}u_{l}).
\]
Thus a similar calculation shows that the corresponding factor in $D_{1}^{\bm{\lambda}_{p}^{(a)}}g_{1}^{\bm{\lambda}_{p}^{(a)}}$ is given by
\[ 
(-1)^{np}q^{\frac{n^{2}(n-1)(e-2)}{2}}\prod_{\stackrel{\scriptstyle 2\leq l\leq e}{l\not= p}}\frac{1}{u_{p}-u_{l}}\prod_{i=1}^{n}\prod_{h=i-n}^{0}(q^{h}u_{p}-u_{l}).
\]

This completes the proof.
\end{proof}
\end{lem}

\begin{lem} We have 
\[
D^{\bm{\lambda}_{p}^{(a)}}g^{\bm{\lambda}_{p}^{(a)}}=(-1)^{a+en}C(e,p)\frac{(1-q)^{n}}{1-q^{n}}\prod_{\stackrel{\scriptstyle 1\leq l\leq e}{l\not= p}}u_{l}^{n}
\prod_{\stackrel{\scriptstyle 1\leq l\leq e}{l\not= p}}\frac{1}{u_{p}-u_{l}}.
\]
\begin{proof}
The proof is divided into five parts. We use Lemma 6.2 and Lemma 6.3.\\
\linebreak
(1) The factors which have the form $q^{*}-q^{*}$. It is given by
\[
\frac{\displaystyle (1-q)^{n}\prod_{\stackrel{\scriptstyle \alpha ,\alpha '\in \bold{A}_{1}}{\alpha >\alpha '}}(q^{\alpha }-q^{\alpha '})
\prod_{\stackrel{\scriptstyle 2\leq l\leq e}{l\not= p}}\prod_{\stackrel{\scriptstyle \alpha ,\alpha '\in \bold{A}_{l}}{\alpha >\alpha '}}(q^{\alpha }-q^{\alpha '})}
{\displaystyle (1-q^{n})\prod_{\alpha \in \bold{A}_{1}}\prod_{h=1}^{\alpha }(q^{h}-1)\prod_{\stackrel{\scriptstyle 2\leq l\leq e}{l\not= p}}\prod_{\alpha \in \bold{A}_{l}}
\prod_{h=1}^{\alpha }(q^{h}-1)}.
\]
By Lemma 6.1 (1) and (2) we have the following two identities:
\begin{center}
$\displaystyle \prod_{\stackrel{\scriptstyle \alpha ,\alpha '\in \bold{A}_{1}}{\alpha >\alpha '}}(q^{\alpha }-q^{\alpha '})
=q^{\frac{n(n+1)(n-1)}{6}}\prod_{\alpha \in \bold{A}_{1}}\prod_{h=1}^{\alpha }(q^{h}-1)$,\\
$\displaystyle \prod_{\stackrel{\scriptstyle \alpha ,\alpha '\in \bold{A}_{l}}{\alpha >\alpha '}}(q^{\alpha }-q^{\alpha '})
=q^{\frac{n(n-1)(n-2)}{6}}\prod_{\alpha \in \bold{A}_{l}}\prod_{h=1}^{\alpha }(q^{h}-1)$.
\end{center}
Thus this part is given by
\[
q^{n(n-1)(en-n-2e+5)}\frac{(1-q)^{n}}{1-q^{n}}.
\]
(2) The factors which have the form $q^{*}u_{k}-q^{*}u_{l}$. We divide into four cases.\\
\linebreak
(2-1) The factors which have the form $q^{*}u_{p}-q^{*}u_{1}$.

This is given by
\[
\frac{\displaystyle \frac{1}{u_{p}-u_{1}}\prod_{i=1}^{n}\prod_{h=i-n}^{1}(q^{h}u_{p}-u_{1})}
{\displaystyle (q^{-1}u_{1}-u_{p})^{n}\prod_{\alpha \in \bold{A}_{1}}\prod_{h=1}^{\alpha }(q^{h-1}u_{1}-u_{p})}.
\]
Now the denominator can be calculated as follows:\\
\linebreak
$\displaystyle (q^{-1}u_{1}-u_{p})^{n}\prod_{i=1}^{n+1}\prod_{h=1}^{n-i+1}(q^{h-1}u_{1}-u_{p})$\\
$\displaystyle =\prod_{i=1}^{n}\prod_{h=0}^{n-i+1}(q^{h-1}u_{1}-u_{p})$\\
$\displaystyle =(-1)^{\sum_{i=1}^{n}(n-i+2)}q^{\sum_{i=1}^{n}(\frac{(n-i)(n-i+1)}{2}-1)}\prod_{i=1}^{n}\prod_{h=0}^{n-i+1}(q^{1-h}u_{p}-u_{1})$\\
$\displaystyle =(-1)^{\frac{n(n+3)}{2}}q^{\frac{n^{3}-7n}{6}}\prod_{i=1}^{n}\prod_{h=i-n}^{1}(q^{h}u_{p}-u_{1}).$\\
\linebreak
Thus this part is given by
\[
(-1)^{\frac{n(n-1)}{2}}q^{-\frac{n^{3}-7n}{6}}\frac{1}{u_{p}-u_{1}}.
\]
(2-2) The factors which have the form $q^{*}u_{p}-q^{*}u_{l}\ (l\geq 2,l\not= p)$. 

This is given by
\[
\frac{\displaystyle \prod_{\stackrel{\scriptstyle 2\leq l\leq e}{l\not= p}}\frac{1}{u_{p}-u_{l}}\prod_{i=1}^{n}\prod_{h=i-n}^{0}(q^{h}u_{p}-u_{l})}
{\displaystyle \prod_{p<l\leq e}(u_{p}-u_{l})^{n}\prod_{2\leq k<p}(u_{k}-u_{p})^{n}\prod_{\stackrel{\scriptstyle 2\leq l\leq e}{l\not= p}}\prod_{\alpha \in \bold{A}_{l}}
\prod_{h=1}^{\alpha }(q^{h}u_{l}-u_{p})}.
\]
Now we can calculate as follow:\\
\linebreak
$\displaystyle \prod_{\alpha \in \bold{A}_{l}}\prod_{h=1}^{\alpha }(q^{h}u_{l}-u_{p})$\\
$\displaystyle =\prod_{i=1}^{n}\prod_{h=1}^{n-i}(q^{h}u_{l}-u_{p})$\\
$\displaystyle =(-1)^{\sum_{i=1}^{n}(n-i)}q^{\sum_{i=1}^{n}\frac{(n-i)(n-i+1)}{2}}\prod_{i=1}^{n}
\prod_{h=1}^{n-i}(q^{-h}u_{p}-u_{l})$\\
$\displaystyle =(-1)^{\frac{n(n-1)}{2}}q^{\frac{n(n+1)(n-1)}{6}}(u_{p}-u_{l})^{n}\prod_{i=1}^{n}\prod_{h=i-n}^{0}(q^{h}u_{p}-u_{l}).$\\
\linebreak
From this calculation one can easily show that this part is given by
\[
(-1)^{\frac{n(n-1)(e-2)}{2}+np}q^{-\frac{n(n+1)(n-1)(e-2)}{6}}\prod_{\stackrel{\scriptstyle 2\leq l\leq e}{l\not= p}}\frac{1}{u_{p}-u_{l}}.
\]
(2-3) The factors which have the form $q^{*}u_{1}-q^{*}u_{l}\ (l\not= 1,p)$. 

This is given by
\[
\frac{\displaystyle \prod_{\stackrel{\scriptstyle 2\leq l\leq e}{l\not= p}}\prod_{(\alpha ,\alpha ')\in \bold{A}_{1}\times \bold{A}_{l}}(q^{\alpha -1}u_{1}-q^{\alpha '}u_{l})}
{\displaystyle\prod_{\stackrel{\scriptstyle 2\leq l\leq e}{l\not= p}} (q^{-1}u_{1}-u_{l})^{n}\prod_{\stackrel{\scriptstyle 2\leq l\leq e}{l\not= p}}\prod_{\alpha \bold{A}_{1}}\prod_{h=1}^{\alpha }(q^{h-1}u_{1}-u_{l})
\prod_{\stackrel{\scriptstyle 2\leq l\leq e}{l\not= p}}\prod_{\alpha \in \bold{A}_{l}}\prod_{h=1}^{\alpha }(q^{h}u_{l}-q^{-1}u_{1})}.
\]
Now one can check easily the following three identities:
\begin{center}
$\displaystyle \prod_{(\alpha ,\alpha ')\in \bold{A}_{1}\times \bold{A}_{l}}(q^{\alpha -1}u_{1}-q^{\alpha '}u_{l})=q^{\frac{n(n+1)(n-2)}{2}}
\prod_{i=1}^{n+1}\prod_{j=1}^{n}(u_{1}-q^{i-j}u_{l})$,\\
$\displaystyle \prod_{\alpha \in \bold{A}_{1}}\prod_{h=1}^{\alpha }(q^{h-1}u_{1}-u_{l})=q^{\frac{n(n+1)(n-1)}{6}}\prod_{i=1}^{n+1}\prod_{h=1}^{n-i+1}(u_{1}-q^{1-h}u_{l})$,\\
$\displaystyle \prod_{\alpha \in \bold{A}_{l}}(q^{h}u_{l}-q^{-1}u_{1})=\frac{\displaystyle (-1)^{\frac{n(n-1)}{2}}q^{-\frac{n(n-1)}{2}}\prod_{i=1}^{n+1}\prod_{h=1}^{i-1}(u_{1}-q^{h+1}u_{l})}{\displaystyle \prod_{h=1}^{n}(u_{1}-q^{h+1}u_{l})}$.
\end{center}
From these calculations we find that this part is given by
\[
(-1)^{\frac{n(n-1)(e-2)}{2}}q^{\frac{n(n+1)(n-1)(e-2)}{3}}.
\]
(2-4) The factors which have the form $q^{*}u_{k}-q^{*}u_{l}\ (k,l\geq 2,\ k\not=l,p)$.

This is given by
\[
\frac{\displaystyle \prod_{\stackrel{\scriptstyle 2\leq k<l\leq e}{k,l\not= p}}\prod_{(\alpha ,\alpha ')\in \bold{A}_{k}\times \bold{A}_{l}}(q^{\alpha }u_{k}-q^{\alpha '}u_{l})}
{\displaystyle\prod_{\stackrel{\scriptstyle 2\leq k<l\leq e}{k,l\not= p}}(u_{k}-u_{l})^{n}\prod_{\stackrel{\scriptstyle k,l\geq 2}{k,l\not= p,\ k\not= l}}\prod_{\alpha \in \bold{A}_{k}}
\prod_{h=1}^{\alpha }(q^{h}u_{k}-u_{l})}.
\]

Let $k,l$ be two integers such that $2\leq k<l\leq e$. We first note the following identity:
\[
\prod_{(\alpha ,\alpha ')\in \bold{A}_{k}\times \bold{A}_{l}}(q^{\alpha }u_{k}-q^{\alpha '}u_{l})=q^{\frac{n^{2}(n-1)}{2}}\prod_{i=1}^{n}\prod_{j=1}^{n}(u_{k}-q^{i-j}u_{l}).
\]
On the other hand, we have:\\
\linebreak
$\displaystyle \prod_{\stackrel{\scriptstyle k,l\geq 2}{k,l\not= p,\ k\not= l}}\prod_{\alpha \in \bold{A}_{k}}
\prod_{h=1}^{\alpha }(q^{h}u_{k}-u_{l})$\\
$\displaystyle =\prod_{\stackrel{\scriptstyle 2\leq k<l\leq e}{k,l\not= p}}\Bigg( \prod_{\alpha \in \bold{A}_{k}}\prod_{h=1}^{\alpha }(q^{h}u_{k}-u_{l})\cdot 
\prod_{\alpha \in \bold{A}_{l}}\prod_{h=1}^{\alpha }(q^{h}u_{l}-u_{k})\Bigg)$\\
$\displaystyle =\prod_{\stackrel{\scriptstyle 2\leq k<l\leq e}{k,l\not= p}}\Bigg( (-1)^{\frac{n(n-1)}{2}}q^{\frac{n(n+1)(n-1)}{6}}\prod_{i=1}^{n}\prod_{h=1}^{n-i}
(u_{k}-q^{-h}u_{l})\prod_{i=1}^{n}\prod_{h=1}^{i-1}(u_{k}-q^{h}u_{l})\Bigg)$.\\
\linebreak
From these computations we find that this part is given by
\[
(-1)^{\frac{n(n-1)}{2}\binom{e-2}{2}}q^{\frac{n(n-1)(2n-1)}{6}\binom{e-2}{2}}.
\]

By combining (2-1)$,\cdots,$(2-4), we can conclude that the factors which have the form $q^{*}u_{k}-q^{*}u_{l}$ is given by
\[
(-1)^{\binom{n}{2}+\binom{n}{2}\binom{e-2}{2}+np}q^{\frac{n(n+1)(n-1)}{6}(e-3)+n+\frac{n(n-1)(2n-1)}{6}\binom{e-2}{2}}\prod_{\stackrel{\scriptstyle
1\leq l\leq e}{l\not =p}}\frac{1}{u_{p}-u_{l}}.
\]
(3) The power of $q$. \\
\linebreak
The exponent is given by
\[
-n-f(n,e)-\binom{n+1}{2}+\sum_{\alpha \in \bold{A}_{1}}\alpha +\frac{n(n-1)(n-2)}{6}+\frac{n(n-1)(en-n+1)}{2}.
\]
A direct computation shows that this is equal to
\[
\frac{n(n-1)(3en-2n+1)}{6}-f(n,e)-n.
\]
(4) The sign.\\
\linebreak
This is given by
\[
(-1)^{a+n+np+\binom{e}{2}\binom{n}{2}+n(e-1)}.
\]
(5) The other factor.\\
\linebreak
This is given by
\[
C(e,p)\times \prod_{k=1}^{e}u_{k}^{n}\times u_{1}^{\binom{n+1}{2}}\times \prod_{\stackrel{\scriptstyle 2\leq l\leq e}{l\not =p}}u_{l}^{n}\times u_{p}^{\binom{n}{2}}
\times u_{1}^{-\sum_{\alpha \in \bold{A}_{1}}}\times \prod_{\stackrel{\scriptstyle 2\leq l\leq e}{l\not =p}}u_{l}^{-\sum_{\alpha \in \bold{A}_{p}}},
\]
and is equal to
\[
C(e,p)\prod_{\stackrel{\scriptstyle 1\leq l\leq e}{l\not =p}}u_{l}^{n}.
\]

Now our desired formula follows from $(1),\cdots, (5)$ and Remark 3.2 directly.

This completes the proof.
\end{proof}
\end{lem}

Next, we must compute $D^{\bm{\lambda}_{1}^{(a)}}g^{\bm{\lambda}_{1}^{(a)}}$. However, the computation of $D^{\bm{\lambda}_{1}^{(a)}}g^{\bm{\lambda}_{1}^{(a)}}$ is quite similar
to (and relatively easier than) that of $D^{\bm{\lambda}_{p}^{(a)}}g^{\bm{\lambda}_{p}^{(a)}}
\ (p\in \{ 2,\ldots,e\})$. So we omit the details. Instead, we write down the corresponding results for reader's convenience.

In this case, $\alpha _{i,j}$ are given by
\begin{eqnarray}
\alpha_{1,i}= \left\{ \begin{array}{ll}
a+n+1 & (i=1) \\
n-i+2 & (2 \leq i\leq b+1) \\
n-i+1 &  (b+2\leq i\leq  n+1),\\
\end{array} \right. \nonumber
\end{eqnarray}
\[
\alpha_{l,i}=n-i\ \ (2\leq l\leq e,1\leq i\leq n).
\]

The following four lemmas correspond to Lemma 6.1, 6.2, 6.3 and 6.4 respectively.

\begin{lem}
For $i,j,l\geq 2$ we have the following.
\begin{itemize}
\item[(1)] $\displaystyle \prod_{j'=1}^{j-1}(q^{\alpha _{l,j'}}-q^{\alpha _{l,j}})=q^{(j-1)(n-j)}\prod_{h=1}^{j-1}(q^{h}-1)$
\item[(2)] $\displaystyle \prod_{i'=1}^{i-1}(q^{\alpha _{1,i'}}-q^{\alpha _{1,i}})$\\
$\displaystyle = \left\{ \begin{array}{ll}
\displaystyle q^{(i-1)(n-i+1)}\frac{q^{a+i}-1}{q^{i-b-1}-1} \prod_{h=1}^{i-1}(q^{h}-1) & (b+3\leq i\leq n+1) \\
\displaystyle q^{(i-1)(n-i+1)}\frac{q^{a+i}-1}{q-1} \prod_{h=1}^{i-1}(q^{h}-1) & (i=b+2) \\
\displaystyle q^{(i-1)(n-i+1)+i-1}\frac{q^{a+i-1}-1}{q^{i-1}-1} \prod_{h=1}^{i-1}(q^{h}-1)  &  (2\leq i\leq b+1).\\
\end{array} \right. \nonumber $
\end{itemize}
\end{lem}
\begin{lem}We have
\[
 D_{0}^{\bm{\lambda}_{p}^{(a)}}g_{0}^{\bm{\lambda}_{p}^{(a)}}=(-1)^{a+n}C(e,1)q^{\frac{n(n+1)(n-1)}{6}}\frac{(1-q)^{n}}{1-q^{n}}.
\]
\end{lem}

\begin{lem}We have
\[
 D_{1}^{\bm{\lambda}_{1}^{(a)}}g_{1}^{\bm{\lambda}_{1}^{(a)}}=q^{\frac{n(n+1)(n-1)(e-1)}{2}}\prod_{2\leq l\leq e}\Bigg( \frac{1}{u_{1}-u_{l}}
\prod_{i=1}^{n}\prod_{h=i-n-1}^{-1}(q^{h}u_{1}-u_{l})\Bigg).
\]
\end{lem}

\begin{lem} We have 
\[
D^{\bm{\lambda}_{1}^{(a)}}g^{\bm{\lambda}_{1}^{(a)}}=(-1)^{a+en}C(e,1)\frac{(1-q)^{n}}{1-q^{n}}\prod_{2\leq l\leq e}\frac{1}{u_{1}-u_{l}}.
\]
\end{lem}

The following is a direct consequence of Lemma 5.3, 6.4 and 6.8.
\begin{thm}
For any mixed braid $\alpha $ with $n$-strands we have\\
\linebreak
$\Delta^{G(e,1)}(\widehat{\alpha })$\\
$\displaystyle =(-1)^{n-1}q^{-\frac{n-wr(\widehat{\alpha })-1}{2}}\frac{1-q}{1-q^{n}}
\sum_{\stackrel{\scriptstyle 0\leq a\leq n-1}{1\leq p\leq e}}(-1)^{a}C(e,p)
\Bigg(\prod_{\stackrel{\scriptstyle 1\leq l\leq e}{l\not =p}}\frac{1}{u_{p}-u_{l}}\Bigg)\chi_{\bm{\lambda}_{p}^{(a)}}(\pi_{n}(\alpha ))$. 
\end{thm}

From this formula we find the relationship between the Alexander polynomial of a mixed link and of the link which is obtained by resolving the twisted parts.

\begin{thm}
Let $\alpha $ be a mixed braid, $\alpha_{0}$ be the braid obtained by avoiding the powers of $t_{0}$ appearing in $\alpha $, and $\Delta (\alpha_{0})$ be the
usual Alexander polynomial of $\alpha_{0}$. Then we have
\[
\Delta^{G(e,1)}(\widehat{\alpha} )=\Delta(\widehat{\alpha_{0}}) \sum_{1\leq p\leq e}C(e,p)u_{p}^{wr_{0}(\alpha )} \Bigg(\prod_{\stackrel{\scriptstyle 1\leq l\leq e}{l\not =p}}
\frac{1}{u_{p}-u_{l}}\Bigg).
\]
Here, $wr_{0}(\alpha )$ is the sum of exponents of $t_{0}$ appearing in $\alpha $.
\begin{proof}
We first remark that by taking $e=1$ we can recover the following formula 
\[
\Delta (\widehat{\alpha_{0}})=(-1)^{n-1}q^{-\frac{n-wr(\widehat{\alpha })-1}{2}}\frac{1-q}{1-q^{n}}
\sum_{ 0\leq a\leq n-1}(-1)^{a} \chi_{\bm{\lambda}_{p}^{(1)}}(\pi_{n}(\alpha_{0}))
\]
of Jones.
Since $\chi _{\bm{\lambda}_{p}^{(a)}}(\alpha )=u_{p}^{wr_{0}(\alpha )}\chi _{\bm{\lambda}_{p}^{(1)}}(\alpha _{0})$,
we have\\
\linebreak
$\Delta^{G(e,1)}(\widehat{\alpha })$\\
$\displaystyle =(-1)^{n-1}q^{-\frac{n-wr(\widehat{\alpha })-1}{2}}\frac{1-q}{1-q^{n}}$\\
$\hfill \displaystyle \times \sum_{\stackrel{\scriptstyle 0\leq a\leq n-1}{1\leq p\leq e}}(-1)^{a}C(e,p)
\Bigg(\prod_{\stackrel{\scriptstyle 1\leq l\leq e}{l\not =p}}\frac{1}{u_{p}-u_{l}}\Bigg)\chi_{\bm{\lambda}_{p}^{(a)}}(\pi_{n}(\alpha ))$\\
$\displaystyle =(-1)^{n-1}q^{-\frac{n-wr(\widehat{\alpha })-1}{2}}\frac{1-q}{1-q^{n}}$\\
$\hfill \displaystyle \times \sum_{1\leq p\leq e}C(e,p)u_{p}^{wr_{0}(\alpha )}
\Bigg(\prod_{\stackrel{\scriptstyle 1\leq l\leq e}{l\not =p}}\frac{1}{u_{p}-u_{l}}\Bigg)
\sum_{0\leq a\leq n-1}(-1)^{a}\chi_{\bm{\lambda}_{p}^{(1)}}(\pi_{n}(\alpha_{0} ))$\\
$\displaystyle =\Delta(\alpha _{0})\sum_{1\leq p\leq e}C(e,p)u_{p}^{wr_{0}(\alpha )} \Bigg(\prod_{\stackrel{\scriptstyle 1\leq l\leq e}{l\not =p}}
\frac{1}{u_{p}-u_{l}}\Bigg).$\\
\linebreak
as desired.
\end{proof}
\end{thm}


\bigskip
\noindent
Graduate~School~of~Science, Osaka~City~University, Sugimoto, Sumiyoshi-ku, Osaka~558-8585, Japan.

\bigskip
\noindent
Email Address: d09saq0L05@ex.media.osaka-cu.ac.jp
\end{document}